\documentclass{birkjour}
 \newtheorem{thm}{Theorem}[section]
 
 \newtheorem{lem}[thm]{Lemma}
 
 \theoremstyle{definition}
 \newtheorem{defn}[thm]{Definition}
 \theoremstyle{remark}
 \newtheorem{rem}[thm]{Remark}
 
 \numberwithin{equation}{section}
 \usepackage[hidelinks]{hyperref}
\usepackage{amssymb,amsmath,amsfonts}
\usepackage{mathtools}
\usepackage{amsthm}
\usepackage{float}
\DeclareMathOperator*{\diag}{diag}

\newcommand{\dd}{\!\mathrm{d}}
\newcommand{\E}{\mathrm{e}}

\usepackage[svgnames]{xcolor}
\usepackage{todonotes}
\usepackage{tikz,pgfplots}
\usetikzlibrary{decorations.markings}
\usetikzlibrary{shapes.geometric}
\pgfdeclarelayer{edgelayer}
\pgfdeclarelayer{nodelayer}
\pgfsetlayers{edgelayer,nodelayer,main}
\tikzstyle{none}=[inner sep=0pt]
\tikzstyle{wh}=[circle,fill=White,draw=Black,line width=0.8 pt]
\tikzstyle{rn}=[circle,fill=Red,draw=Black,line width=0.8 pt]
\tikzstyle{gn}=[circle,fill=Lime,draw=Black,line width=0.8 pt]
\tikzstyle{yn}=[circle,fill=Yellow,draw=Black,line width=0.8 pt]
\tikzstyle{simple}=[-,draw=Black,line width=2.000]
\tikzstyle{arrow}=[-,draw=Black,postaction={decorate},decoration={markings,mark=at position .5 with {\arrow{>}}},line width=2.000]
\tikzstyle{tick}=[-,draw=Black,postaction={decorate},decoration={markings,mark=at position .5 with {\draw (0,-0.1) -- (0,0.1);}},line width=2.000]

\usepackage{subfig}
\usepackage{booktabs}
\hypersetup{breaklinks=true}
\usepackage{xparse}
\usepackage{ifmtarg}
\usepackage{makecell}
\makeatletter

\DeclareDocumentCommand{\s}{o o o}{
f
        \IfNoValueTF{#1}
            {}%
            {\@ifmtarg{#1}{}{^{#1}}}
        \IfNoValueTF{#2}
            {}%
            {\@ifmtarg{#2}{}{_{#2}}}%
}

\DeclareDocumentCommand{\es}{o o o}{
\mathfrak{f}
        \IfNoValueTF{#1}
            {}%
            {\@ifmtarg{#1}{}{^{#1}}}
        \IfNoValueTF{#2}
            {}%
            {\@ifmtarg{#2}{}{_{#2}}}%
}

\begin{document}

%-------------------------------------------------------------------------
% editorial commands: to be inserted by the editorial office
%
%\firstpage{1} \volume{228} \Copyrightyear{2004} \DOI{003-0001}
%\seriesextraline{This is the Concrete Title of this Book\br H.E. R and S.T.C. W, Eds.}
%\firstpage{1}

\title[Spectral analysis of matrix-sequences via GLT momentary symbols]
 {A note on the spectral analysis of matrix-sequences via GLT momentary symbols: from all-at-once solution of parabolic problems to distributed fractional order matrices}

%----------Author 1
\author[M. Bolten]{Matthias Bolten}

\address{%
Bergische Universit\"{a}t Wuppertal\\
Department of Mathematics and Computer Science\\
 Gau\ss stra\ss e 20, 42119 Wuppertal,\\ Germany}

\email{bolten@math.uni-wuppertal.de}

%----------Author 2
\author[S.-E. Ekstr\"{o}m]{Sven-Erik Ekstr\"{o}m}
\address{Uppsala University\\ Department of Information Technology\\  Lägerhyddsv. 2, hus 2, SE-751 05, Uppsala\\ Sweden}
\email{sven-erik.ekstrom@it.uu.se}

%----------Author 3
\author[I. Furci]{Isabella Furci}

\address{%
Bergische Universit\"{a}t Wuppertal\\
Department of Mathematics and Computer Science\\
 Gau\ss stra\ss e 20, 42119, Wuppertal\\ Germany}

\email{furci@uni-wuppertal.de}

%----------Author 4
\author[S. Serra-Capizzano]{Stefano Serra-Capizzano}
\address{Insubria University\\ Department of Science and High Technology\\
 via Valleggio 11, 22100, Como\\ Italy.}
\email{s.serracapizzano@uninsubria.it}
\address{Uppsala University\\ Department of Information Technology\\  Lägerhyddsv. 2, hus 2, SE-751 05, Uppsala\\ Sweden}
\email{stefano.serra@it.uu.se}




%%% ----------------------------------------------------------------------

\begin{abstract}
The first focus of this paper is the characterization of the spectrum and the singular values of the coefficient matrix stemming from the discretization with space-time grid for a parabolic diffusion problem and from the approximation of distributed order fractional equations.
For this purpose we will use the classical GLT theory and the new concept of GLT momentary symbols. The first permits to describe the
singular value or eigenvalue asymptotic distribution of the sequence of the coefficient matrices, the latter permits to derive a function, which describes the singular value or eigenvalue distribution of the matrix of the sequence, even for small matrix-sizes but under given assumptions. The note is concluded with a list of open problems, including the use of our machinery in the study of iteration matrices, especially those concerning multigrid-type techniques.
\end{abstract}

\maketitle

\section{Introduction and notation}
\label{sec:introduction}

As well known, many practical applications require to solve numerically linear systems of Toeplitz kind and of large dimensions. As a consequence a number of  iterative techniques, such as preconditioned Krylov methods, multigrid procedures, and sophisticated combination of them have been designed (see \cite{ChJi2007,Ng2004} and the references therein). Linear systems with Toeplitz coefficient matrices of large dimension arise when dealing with the numerical solution of (integro-)differential equations and of problems with Markov chains. More recently, new examples of real world problems have emerged.
\textcolor{black}{
The first focus of this paper is the characterization of the spectrum and the singular values of the coefficient matrix stemming from the discretization with a space-time grid for a parabolic diffusion problem.
More specifically, we consider the diffusion equation in one space dimension,
\begin{equation*}
u_t=u_{xx}, \quad x\in (a,b),\ t\in [0, T],
\end{equation*}
and we approximate our parabolic model problem on a rectangular space-time grid consisting of $N_t$ time intervals and $N_x$ space intervals. 
}

The second focus  concerns  the matrix-sequences involved with the  discretization of distributed order fractional differential equations (FDEs) which have gained a lot of attention. Owing to the nonlocal nature of fractional operators, independently of the locality of the approximation methods, the matrix structures are dense and under assumptions of uniform step-sizing and of constant coefficients in the involved operators, the matrices are again of Toeplitz type (unilevel, or multilevel according to the dimensionality of the considered domains).

When the fractional order is fixed, the spectral analysis of such matrices (conditioning, extremal eigenvalues etc) can be performed, by exploiting the well-established analysis of the spectral features of Toeplitz matrix-sequences generated by Lebesgue integrable functions and the more recent Generalized Locally Toeplitz (GLT) theory \cite{GLT-bookI}; see for instance \cite{DoMa2016, DoMa2018}. However in the case of the numerical approximation of distributed-order fractional operators, also the spectral analysis of the resulting matrices is more involved. We recall that distributed-order FDEs can be interpreted as a parallel distribution of derivatives of fractional orders, whose most immediate application consists in the physical modeling of systems characterized by a superposition of different processes operating in parallel. As an example, we mention the application of fractional distributed-order operators as a tool for accounting memory effects in composite materials \cite{CaFa2017} or multi-scale effects \cite{CaGi2018}. For a detailed review on the topic we refer the reader to \cite{DiPa2021}.

\textcolor{black}{In order to study the involved structured linear systems  of both integral and differential equations, we will use the classical theory of GLT matrix-sequences \cite{GLT-bookI,GLT-bookII} and the new concept of GLT momentary symbols.} The first permits to describe the singular value or eigenvalue asymptotic distribution of the sequence of the coefficient matrices, the latter permits to derive a function, which describes the singular value or eigenvalue distribution of a fixed matrix of the sequence, even for small matrix-sizes, but under given assumptions.

This paper is organized as follows. The remaining part of this section is devoted to definitions, notation, and to the necessary background for our analysis: in particular we provide a formal definition of  GLT momentary symbols. Section \ref{sec:problem} is devoted to setting up the problem and to derive the relevant matrix structures. The distributional analysis both for the eigenvalues and singular values is the main focus of Subsection \ref{sec:coefficientmatrix}, while Section \ref{sec:fractional} contains similar results for specific matrix-structures with generating function depending on the matrix-size and which arise in the context of fractional differential equations with distributed orders.
Section \ref{sec:conclusions}  contains conclusions and a list of open problems, including the use of our machinery in the study of iteration matrices, especially those concerning multigrid-type techniques.

\subsection{Background and definitions}
\label{sec:introduction:background}
Throughout this paper, we will use the following notations. Let ${f}:G\to\mathbb{C}$ be a function belonging to $L^1(G)$, with $G\subseteq\mathbb R^\ell$, $\ell\ge 1$, a measurable set.
We denote by $\{A_{n}\}_{n}$ the matrix-sequence, whose elements are given by the matrices $A_{n}$ of dimension $n \times n$. Let $s,d \in \mathbb{N}$. Let $\mathbf{n}=(n_1,n_2,\dots,n_d)$ be a multi-index, we indicate by  $\{A_{\mathbf{n}}\}_{\mathbf{n}}$, the $d$-level $s\times s$ block matrix-sequence, whose elements are the matrices $A_\mathbf{n}$ of size $d=d(\mathbf{n},s)=sn_1n_2\cdots n_d$.
\subsection{Toeplitz and circulant matrix-sequences}
\label{sec:introduction:tep_circ}
In the following we report the main background concerning the concepts of Toeplitz and circulant matrices, for simplicity, in the scalar unilevel setting. We only provide the generalization in the block multilevel case of the results that will be exploited for the purpose of the paper.
\begin{defn}
\label{def:toeplitz_scalar}
An $n\times n$ Toeplitz matrix $A_n$ is a matrix that has equal {entries} along each diagonal, and can be written as
\begin{equation*}
A_n=\left[a_{i-j}\right]_{i,j=1}^{n}=\left[\begin{smallmatrix}
a_0 & a_{-1} & a_{-2} & \cdots & \cdots & a_{1-n}\vphantom{\ddots}\\
a_1 & \ddots & \ddots & \ddots & & \vdots\\
a_2 & \ddots & \ddots & \ddots & \ddots & \vdots\\
\vdots & \ddots & \ddots & \ddots & \ddots & a_{-2}\\
\vdots & & \ddots & \ddots & \ddots & a_{-1}\\
a_{n-1} & \cdots & \cdots & a_2 & a_1 & a_0\vphantom{\ddots}\\
&
\end{smallmatrix}\right], \ \ \ a_j\in \mathbb{C}, \ j=1-n,\ldots,n-1.
\end{equation*}
\end{defn}
In the following we focus on the two important sub-classes given by the Toeplitz matrices $T_{n}(f) \in \mathbb{C}^{n \times n}$ and the circulant matrices $C_{n}(f) \in \mathbb{C}^{n \times n}$,  associated with a function $f$, called the \textbf{generating function.}

\begin{defn}\label{def:toeplitz_scalar_generating}
Given $f$ belonging to $L^1([-\pi,\pi])$ and periodically extended to the whole real line, the matrix $T_{n}(f)$ is defined as
\begin{equation*}
  T_n(f)=\left[\hat f_{i-j}\right]_{i,j=1}^n,
\end{equation*}
where
\begin{equation}
  \hat{f}_{k}\coloneqq\frac1{2\pi}\int_{-\pi}^{\pi}\!\!f(\theta)\,\E^{-k\mathbf{i} \theta}\dd\theta,\qquad k\in\mathbb Z,\qquad \mathbf{i}^2=-1,\label{eq:introduction:background:fourier}
\end{equation}
are the Fourier coefficients of $f$, and
\begin{equation*}
  f(\theta)=\!\!\sum_{k=-\infty}^{\infty}\!\!\hat{f}_{k}\E^{k\mathbf{i} \theta},\label{eq:introduction:fourierseries}
\end{equation*}
is the Fourier series of $f$.
\end{defn}

{\begin{defn}\label{def:Cir}
Let the Fourier coefficients of a given function ${f}\in L^1([-\pi,\pi])$ be defined as in formula (\ref{eq:introduction:background:fourier}).
Then, we can define the $n \times n$
circulant matrix $C_{n}(f)$ associated with $f$, as
 \begin{equation}
 C_{n}(f)=\!\!\!\!\!\!\sum_{j=-(n-1)}^{n-1}\!\!\!\!\!\!\hat{a}_{j}Z_{n}^{j}=\mathbb{F}_{n}  D_{n}(f) \mathbb{F}_{n}^{*},\label{eq:introduction:background:circulant:schur}
 \end{equation} where $^*$ denotes the transpose conjugate,  $Z_{n}$ is the $n \times n$ matrix defined by
\begin{equation*}
\left(Z_{n}\right)_{ij}=\begin{cases}
1,&\text{if }\mathrm{mod}(i-j,n)=1,\\
0,&\text{otherwise}.
\end{cases}
\end{equation*}
Moreover,
\begin{equation*}
  D_{n}(f)=\diag\left(s_n(f(\theta_{j,n}^c))\right),\quad j=1,\ldots,n,\label{eig-circ}
\end{equation*}
where
\begin{equation}
   \theta_{j,n}^c=\frac{(j-1)2\pi}{n},\quad j=1,\ldots,n,\label{eq:introduction:background:circulant:grid-circ}
 \end{equation}
and $s_{n}(f(\theta))$ is the $ n$th Fourier sum of $f$ given by
\begin{equation*}
s_{n}(f({\theta}))= \sum_{k=1-n}^{n-1}  \hat{f}_{k}
\E^{k\mathbf{i}\theta}.\label{fourier-sum}
\end{equation*}
The matrix $\mathbb{F}_n$ is the so called Fourier matrix of order $n$, given by
 \begin{align*}
 (\mathbb{F}_{n})_{i,j}=\frac{1}{\sqrt{n}} \E^{\mathbf{i}(i-1)\theta_{j,n}^c}, \quad i,j=1,\ldots,n.
 \end{align*}
 In the case of the Fourier matrix, we have $\mathbb{F}_{n}\mathbb{F}_{n}^*=\mathbb{I}_{n}$, that is $\mathbb{F}_{n}$ is complex-symmetric and unitary, with $\mathbb{I}_{n}$ being the identity of size $n$.

The proof of the second equality in (\ref{eq:introduction:background:circulant:schur}), which implies that the columns of the Fourier matrix $\mathbb{F}_n$  are the eigenvectors of $C_{n}(f)$, can be found in \cite[Theorem 6.4]{GLT-bookI}.

\end{defn}
Note that, from the definition follows that if ${f}$ is a trigonometric polynomial of fixed degree less than $n$,   the entries of $D_n(f)$ are the eigenvalues of $C_n(f)$, explicitly given by sampling the generating function $f$ {using} the grid $\theta_{j,n}^c$.
 \begin{equation}
 \begin{split}
  \lambda_j(C_n(f))&=f\left(\theta_{j,n}^c\right),\quad j=1,\ldots,n,\nonumber\\
   D_{n}(f)&=\diag\left(f\left(\theta_{j,n}^c\right)\right),\quad j=1,\ldots,n. \label{eq:introduction:background:circulant:eig-circSEE}
   \end{split}
 \end{equation}

 The type of domain (either one-dimensional $[-\pi,\pi]$ or d-dimensional $[-\pi,\pi]^d$) and codomain (either the complex field or the space of $s \times s$ complex matrices) of $f$ gives rise to different kinds of Toeplitz matrices, see Table~\ref{ttt} for a complete overview.

\begin{table}[]
	\begin{center}
	\small
  \caption{Different types of generating function and the associated Toeplitz matrix.}
  \label{ttt}
		\begin{tabular}{ll|ll}
			\toprule
			\multicolumn{2}{c|}{\textbf{Type of generating function}} & \multicolumn{2}{c}{\textbf{Associated Toeplitz matrix}}\\
			\midrule
			univariate scalar & $f(\theta):[-\pi,\pi]\to \mathbb{C}$ & unilevel scalar & $T_{n}(f)\in \mathbb{C}^{n\times n}$\\
			$d$-variate scalar & $f(\boldsymbol{\theta}):[-\pi,\pi]^d\to \mathbb{C}$ & $d$-level scalar & $T_{\mathbf{n}}(f)\in \mathbb{C}^{d(\mathbf{n},1)\times d(\mathbf{n},1)}$\\
			univariate matrix-valued & $\mathbf{f}(\theta):[-\pi,\pi]\to \mathbb{C}^{s\times s}$& unilevel block & $T_{n}(\mathbf{f})\in \mathbb{C}^{d({n},s)\times d({n},s)}$\\
			$d$-variate matrix-valued & $\mathbf{f}(\boldsymbol{\theta}):[-\pi,\pi]^d\to \mathbb{C}^{s\times s}$ & $d$-level block & $T_{\mathbf{n}}(\mathbf{f})\in \mathbb{C}^{d(\mathbf{n},s)\times d(\mathbf{n},s)}$\\
			\bottomrule
		\end{tabular}
	\end{center}
	
\end{table}

In particular, we provide the definition of a $d$-level $s\times s$ block Toeplitz matrices $T_\mathbf{n}(\mathbf{f})$ starting from  $d$-variate matrix-valued function $\mathbf{f}:[-\pi,\pi]^{d}\rightarrow \mathbb{C}^{s\times s}$, with
with $\mathbf{f}\in L^1([-\pi,\pi]^d)$.
\begin{defn}\label{def:toeplitz_block_generating}
Given a function $\mathbf{f}:[-\pi,\pi]^{d}\rightarrow \mathbb{C}^{s\times s}$ its Fourier coefficients are given by
\begin{equation*}
  \hat{\mathbf{f}}_{\mathbf{k}}\coloneqq
  \frac1{(2\pi)^d}
  \int_{[-\pi,\pi]^d}\mathbf{f}(\boldsymbol{\theta})\E^{-\mathbf{i}\left\langle {\mathbf{k}},\boldsymbol{\theta}\right\rangle}\mathrm{d}\boldsymbol{\theta}\in\mathbb{C}^{s\times s},
  \qquad \mathbf{k}=(k_1,\ldots,k_d)\in\mathbb{Z}^d,\label{fhat}
\end{equation*}
where $\boldsymbol{\theta}=(\theta_1,\ldots,\theta_d)$, $\left\langle \mathbf{k},\boldsymbol{\theta}\right\rangle=\sum_{i=1}^dk_i\theta_i$, and the integrals of matrices are computed elementwise. The associated generating function an be defined via its Fourier series as
\begin{equation*}
\mathbf{f}(\boldsymbol{\theta})=\sum_{\mathbf{k}\in \mathbb{Z}^d}\hat{\mathbf{f}}_{\mathbf{k}}\E^{\mathbf{i}\left\langle {\mathbf{k}},\boldsymbol{\theta}\right\rangle}.\label{eq:introduction:matrixvaluedsymbol}
\end{equation*}

The $d$-level $s\times s$ block Toeplitz matrix associated with $\mathbf{f}$ is the matrix of dimension $d({\mathbf{n},s})$, where $\mathbf{n}=(n_1,\ldots,n_d)$,  given by
\begin{equation*}
T_\mathbf{n}(\mathbf{f})=
\sum_{\mathbf{e}-\mathbf{n}\le \mathbf{k}\le \mathbf{n}-\mathbf{e}} T_{n_1}(\E^{\mathbf{i}k_1\theta_1})\otimes \cdots \otimes
T_{n_d}(\E^{\mathbf{i}k_d\theta_1})\otimes \hat{\mathbf{f}}_{\mathbf{k}},
\end{equation*}
where $\mathbf{e}$ is the vector of all ones and where $\mathbf{s}\le \mathbf{t}$ means that $s_j\le t_j$ for any $j=1,\ldots,d$.
\end{defn}

\begin{defn}
If $\textbf{n}\in \mathbb{N}^d$ and $\textbf{a} : [0,1]^d \to \mathbb{C}^{s\times s}$, we define the $\textbf{n}$-th $d$-level and $s\times s$ block diagonal sampling matrix as the following multilevel block diagonal matrix of dimension $d(\textbf{n},s)$:
\begin{equation*}\label{eq:def_diagonal_sampling}
D_{\textbf{n}}(\textbf{a}) = \diag_{\textbf{e}\le \textbf{j} \le \textbf{n}} \textbf{a}\left(\frac{\textbf{j}}{\textbf{n}}\right),
\end{equation*}
where we recall that $\textbf{e}\le \textbf{j} \le \textbf{n}$ means that $\textbf{j}$ varies from $\textbf{e}$ to $\textbf{n}$ following the lexicographic ordering.

\end{defn}

The following result provides an important relation between tensor products and multilevel Toeplitz matrices.
\begin{lem}{\rm \cite{GLT-bookII}}\label{lemm:tensor_prod}
Let $f_1 , \dots, f_d \in L^1([-\pi,\pi])$,  $\mathbf{n}=(n_1,n_2, \dots, n_d) \in \mathbb{N}^d$. Then,
\begin{equation*}
 T_{ n_1} ( f_1 ) \otimes \dots \otimes T_{n_d} ( f_d )= T_{\mathbf{n}} ( f_1 \otimes \dots \otimes f_d ),
\end{equation*}
where the Fourier coefficients of $f_1 \otimes \dots \otimes f_d$ are given by
\begin{equation*}
( f_1 \otimes \dots \otimes f_d)_{\mathbf{k}} = ( f_1 )_{k_1} \dots ( f_d )_{k_d}, \quad \mathbf{k} \in \mathbb{Z}^d.
\end{equation*}
\end{lem}

\subsection{Asymptotic distributions}
\label{sec:introduction:distributions}
In this subsection we introduce the definition of \textit{asymptotic distribution} in the sense of the eigenvalues and of the singular values, first for a generic matrix-sequence $\{A_{n}\}_{n}$, and then we report specific results concerning the distributions of Toeplitz and circulant matrix-sequences. Finally, we recall the notion of GLT algebra and we introduce a general notion of GLT momentary symbols. We remind that in a more specific and limited setting, the notion of momentary symbols is given in \cite{Momentary_1}: here we generalize the definition in \cite{Momentary_1}.

\begin{defn} {\rm \cite{GLT-bookI,GLT-bookII,gsz,TyZ}}
\label{def:introduction:background:distribution}
Let $f,{\mathfrak{f}}:G\to\mathbb{C}$ be  measurable functions, defined on a measurable set $G\subset\mathbb{R}^\ell$ with $\ell\ge 1$, $0<\mu_\ell(G)<\infty$.
Let $\mathcal{C}_0(\mathbb{K})$ be the set of continuous functions with compact support over $\mathbb{K}\in \{\mathbb{C}, \mathbb{R}_0^+\}$ and let $\{A_{\mathbf{n}}\}_{\mathbf{n}}$, be a sequence of matrices with eigenvalues $\lambda_j(A_{\mathbf{n}})$, $j=1,\ldots,{d_\mathbf{n}}$ and singular values $\sigma_j(A_{\mathbf{n}})$, $j=1,\ldots,{d_\mathbf{n}}$.
Then,

\begin{itemize}
  \item The matrix-sequence  $\{A_{\mathbf{n}}\}_\mathbf{n}$ is \textit{distributed as $\s$ in the sense of the \textbf{singular values}}, and we write
    \begin{align}
      \{A_{\mathbf{n}}\}_{\mathbf{n}}\sim_\sigma \s,\nonumber
    \end{align}
    if the following limit relation holds for all $F\in\mathcal{C}_0(\mathbb{R}_0^+)$:
		\begin{equation}
		  \lim_{\mathbf{n}\to\infty}\frac{1}{{d_\mathbf{n}}}\sum_{j=1}^{{d_\mathbf{n}}}F(\sigma_j(A_\mathbf{n}))=
		  \frac1{\mu_\ell(G)}\int_G  F({|\s(\boldsymbol{\theta})|})\,\dd{\boldsymbol{\theta}}.\label{eq:introduction:background:distribution:sv}
		\end{equation}
    The function $\s$ is called the \textbf{singular value symbol} which describes asymptotically the singular value distribution of the matrix-sequence $ \{A_{\mathbf{n}}\}_{\mathbf{n}}$.
	\item The matrix-sequence $\{A_{\mathbf{n}}\}_{\mathbf{n}}$ is \textit{distributed as $ \es$ in the sense of the \textbf{eigenvalues}}, and we write
    \begin{equation*}
           \{A_{\mathbf{n}}\}_{\mathbf{n}}\sim_\lambda \es,
    \end{equation*}
    if the following limit relation holds for all $F\in\mathcal{C}_0(\mathbb{C})$:
    \begin{equation}
       \lim_{\mathbf{n}\to\infty}\frac{1}{{d_\mathbf{n}}}\sum_{j=1}^{{d_\mathbf{n}}}F(\lambda_j(A_{\mathbf{n}}))=
      \frac1{\mu_\ell(G)}\int_G \displaystyle  F({\es(\boldsymbol{\theta})})\,\dd{\boldsymbol{\theta}}.\label{eq:introduction:background:distribution:ev}
    \end{equation}
    The function $\mathfrak{f}$ is called the \textbf{eigenvalue symbol} which describes asymptotically the eigenvalue distribution of the matrix-sequence $ \{A_{\mathbf{n}}\}_{\mathbf{n}}$.
  \end{itemize}

\end{defn}

\begin{rem}
\label{rem:introduction:background:1}
 Note that, if $A_\mathbf{n}$ is normal for any $\mathbf{n}$ or at least definitely, then $\{A_{\mathbf{n}}\}_{\mathbf{n}}\sim_\sigma \s$ and  $\{A_{\mathbf{n}}\}_{\mathbf{n}}\sim_\lambda \es$ imply that $\s=\es$.
 Of course this is true for Hermitian Toeplitz matrix-sequences as emphasized in Theorem \ref{thm:prel_toepl_szeg} and Theorem \ref{thm:prel_toepl_szeg_multi}.

 Moreover, considering the case $d=1$, if $\s$ (or $\es$) is smooth enough, then the informal interpretation of the limit relation \eqref{eq:introduction:background:distribution:sv} (or \eqref{eq:introduction:background:distribution:ev}) is that, for $n$ sufficiently large, the $n$ singular values (or eigenvalues) of $A_{n}$ can be approximated by a sampling of $|\s(\theta)|$ (or $\es(\theta)$) on an equispaced grid of the interval~$G$, up to the presence of possibly $o(n)$ outliers. It is worthy to notice that in most of the Toeplitz and PDE/FDE applications the number of actual outliers is often limited to $O(1)$ with often a small number of outliers (see \cite{GLT-bookI,GLT-bookII,GLT-bookIII,GLT-bookIV} and references therein).

 The generalization of Definition \ref{def:introduction:background:distribution} and Remark \ref{rem:introduction:background:1} to the block setting and multilevel block setting can be found in \cite{GLT-bookIII,GLT-bookIV} and in the references therein.
\end{rem}
In the case where the matrix-sequence is a Toeplitz matrix-sequence generated by a function, the singular value distribution and the spectral distribution have been well studied in the past few decades. At the beginning Szeg{\H{o}} in \cite{gsz} showed that the eigenvalues of the Toeplitz matrix $T_n(f)$ generated by real-valued $f\in L^{\infty}([-\pi,\pi])$ are asymptotically distributed as~$\s$. Moreover, under the same assumption on $f$, Avram and Parter \cite{MR952991,MR851935} proved that the singular values of $T_n(f)$ are distributed as $|\s|$. This result has been undergone many generalizations and extensions among the years (see \cite{GLT-bookI,GLT-bookII,GLT-bookIII,GLT-bookIV} and the references therein).

The generalized Szeg{\H{o}} theorem that describes the singular value and spectral distribution of Toeplitz sequences generated by  a scalar $f\in L^1([-\pi, \pi])$ is given as follows \cite{MR1481397}.
\begin{thm}\label{thm:prel_toepl_szeg}
	Suppose $f \in L^{1}([-\pi,\pi])$. Let $T_n(f)$ be the Toeplitz matrix generated by $f$. We have
	\begin{equation*}
	\{T_n(f)\}_n \sim_{\sigma} f.
	\end{equation*}
	Moreover, if $f$ is real-valued almost everywhere (a.e.), then
	\begin{equation*}
	\{T_n(f)\}_n \sim_{\lambda} f.
	\end{equation*}
\end{thm}

Tilli \cite{MR1671591} generalized the proof to the block-Toeplitz setting and we report the extension of the eigenvalue result to the case of multivariate Hermitian matrix-valued generating functions.

\begin{thm}\label{thm:prel_toepl_szeg_multi}
Suppose $\mathbf{f}\in L^1([-\pi,\pi]^d,s)$ with positive integers $d,s$. Let $T_{\bf n}(\mathbf{f})$ be the Toeplitz matrix generated by $\mathbf{f}$. We have
\[
\{T_{\bf n}(\mathbf{f})\}_{{\bf n}}\sim_\sigma~\mathbf{f}.
\]
Moreover, if $\mathbf{f}$ is a Hermitian matrix-valued function a.e., then,
\[
\{T_{\bf n}(\mathbf{f})\}_{{\bf n}}\sim_\lambda~\mathbf{f}.
\]
\end{thm}
Concerning the circulant matrix-sequences, though the eigenvalues of a $C_{\bf n}(\textbf{f})$ are explicitly known, a result like Theorem \ref{thm:prel_toepl_szeg} and Theorem  \ref{thm:prel_toepl_szeg_multi}
does not hold for sequences $\left\{C_{\bf n}(\textbf{f})\right\}_{{\bf n}}$ in general. Indeed, the Fourier sum of $\textbf{f}$ converges to $\textbf{f}$ under quite restrictive assumptions (see \cite{zygmund}). In particular, if $\textbf{f}$ belongs to the Dini-Lipschitz class,
then $\left\{C_{\bf n}(\textbf{f})\right\}_{{\bf n}}\sim_\lambda\textbf{f}$, (see \cite{estatico-serra} for more relationships between circulant sequences and spectral distribution results).

\subsection{Matrix algebras}\label{sec:matrix_algebra}

A part from the circulant algebra, introduced in Section \ref{sec:introduction:tep_circ}, we recall that other particular matrix algebras have interesting properties and can be exploited for our purpose. In particular, we mention  the well-known $\tau$-algebras, see \cite{bozzo} and references therein.
Here,  we restrict the analysis to the case of the
matrix algebras $ {\tau_{\varepsilon,\varphi}}$, introduced in~\cite{bozzo}, where an element of the algebra is a matrix
\begin{align}
T_{n,\epsilon,\varphi}(g)=\left[
\begin{array}{ccccc}
a+\varepsilon b &b\\
b&a&b\\
&\ddots&\ddots&\ddots\\
&&b&a&b\\
&&&b&a+\varphi b\
\end{array}
\right],\nonumber
\end{align}
 We can associate to this matrix a function $g$ of the form $g(\theta)=a+2b\cos\theta$.
For some values of $\varepsilon$ and $\varphi$ the exact eigenvalues of $T_{n,\epsilon,\varphi}(g)$ are given by sampling with specific grids; for detailed examples see~\cite{Momentary_1} and \cite{taueconomy} for asymptotic results.

In Table~\ref{tbl:taugrids} we provide the proper grids $\theta^{(\varepsilon,\varphi)}_{j,n}$ and $\Theta_{i,j,n}^{(\varepsilon,\varphi)}$ to give the exact eigenvalues and eigenvectors respectively for $\varepsilon, \varphi \in\{-1,0,1\}$.
\begin{table}[!ht]
\centering
\caption{Grids for  {$\tau_{\varepsilon,\varphi}$}-algebras,  $\varepsilon,\varphi\in\{-1,0,1\}$; $\theta_{j,n}^{(\varepsilon,\varphi)}$ and $\Theta_{j,n}^{(\varepsilon,\varphi)}$ are the grids used to compute the eigenvalues and eigenvectors, respectively. {The standard naming convention (\texttt{dst-*} and \texttt{dct-*}) in parenthesis; see, e.g., \cite[Appendix 1]{CeccheriniSilberstein2008}.}}
\label{tbl:taugrids}
\begin{tabular}{|r|ccc||c|}
\hline
&&$\theta_{j,n}^{(\varepsilon,\varphi)}$&&$\Theta_{i,j,n}^{(\varepsilon,\varphi)}$\\[0.2em]
\hline
\diaghead{\theadfont MMMMM}{ $\varepsilon$ }{ $\varphi$ }&
\thead{-1}&\thead{0}&\thead{1}&\thead{-1, 0, 1}\\[0.2em]
\hline
\thead{-1}&\makecell{{{\tiny (\texttt{dst-2})}}\\$\frac{j\pi}{n}$}&\makecell{{{\tiny (\texttt{dst-6})}}\\$\frac{j\pi}{n+1/2}$}&\makecell{{{\tiny (\texttt{dst-4})}}\\$\frac{(j-1/2)\pi}{n}$}&$(i-1/2)\theta_{j,n}^{(\varepsilon,\varphi)}$\\[0.2em]
\thead{0}&\makecell{{{\tiny (\texttt{dst-5})}}\\$\frac{j\pi}{n+1/2}$}&\makecell{{{\tiny (\texttt{dst-1})}}\\$\frac{j\pi}{n+1}$}&\makecell{{{\tiny (\texttt{dst-7})}}\\$\frac{(j-1/2)\pi}{n+1/2}$}&$i\theta_{j,n}^{(\varepsilon,\varphi)}$\\[0.2em]
\thead{1}&\makecell{{{\tiny (\texttt{dct-4})}}\\$\frac{(j-1/2)\pi}{n}$}&\makecell{{{\tiny (\texttt{dct-8})}}\\$\frac{(j-1/2)\pi}{n+1/2}$}&\makecell{{{\tiny (\texttt{dct-2})}}\\$\frac{(j-1)\pi}{n}$}&$(i-1/2)\theta_{j,n}^{(\varepsilon,\varphi)}+\frac{\pi}{2}$\\[0.2em]
\hline
\end{tabular}
\end{table}

Since all grids $\theta_{j,n}^{(\varepsilon,\varphi)}$ associated with  {$\tau_{\varepsilon,\varphi}$}-algebras where $\varepsilon,\varphi\in\{-1,0,1\}$ are uniformly spaced grids, we know that
\begin{alignat*}{7}
\theta_{j,n}^{(1,1)}&<\theta_{j,n}^{(0,1)}&&=\theta_{j,n}^{(1,0)}\nonumber\\
&&&<\theta_{j,n}^{(-1,1)}&&=\theta_{j,n}^{(1,-1)}\nonumber\\
&&&&&<\theta_{j,n}^{(0,0)}\nonumber\\
&&&&&<\theta_{j,n}^{(-1,0)}&&=\theta_{j,n}^{(0,-1)}\nonumber\\
&&&&&&&<\theta_{j,n}^{(-1,-1)},\qquad\qquad \forall j=1,\ldots,n.\label{eq:gridcomparison}
\end{alignat*}

\subsection{Theory of Generalized Locally Toeplitz (GLT) sequences}
\label{sec:introduction:glt}
In this subsection we will introduce the main properties from the theory of Generalized Locally Toeplitz (GLT) sequences and the practical features, which are sufficient for our purposes, see~\cite{GLT-bookI,GLT-bookII,GLT-bookIII,GLT-bookIV}.

In particular, we consider the multilevel and block setting with $d$ being the number of levels.
\begin{description}
  \item[GLT1] Each GLT sequence has a singular value symbol $\mathbf{f}(\boldsymbol{\theta,x})$ which is measurable according to the Lebesgue measure and according to the second item in Definition \ref{def:introduction:background:distribution} with $\ell=2d$.  In addition, if the sequence is Hermitian, then the distribution also holds in the eigenvalue sense.

 We specify that a GLT sequence $\{A_\mathbf{n}\}_\mathbf{n}$ has GLT symbol $\mathbf{f}(\boldsymbol{\theta,x})$ writing $\{A_\mathbf{n}\}_\mathbf{n}\sim_{\textsc{glt}} \mathbf{f}(\boldsymbol{\theta,x})$, $(\boldsymbol{\theta,x})\in [-\pi,\pi]^d\times [0,1]^d$.

  \item[GLT2] The set of GLT sequences form a $*$-algebra, i.e., it is closed under linear combinations, products, inversion (whenever the symbol is singular, at most, in a set of zero Lebesgue measure), and conjugation.
  Hence, we obtain the GLT symbol of algebraic operations  of a finite set of GLT sequences by performing the same algebraic manipulations of the symbols of the considered GLT sequences.
	
  \item[GLT3] Every Toeplitz sequence $\{T_\mathbf{n}(\mathbf{f})\}_n$ generated by a function $\mathbf{f}(\boldsymbol{\theta})$ belonging to $L^1([-\pi,\pi]^d)$ is a GLT sequence and with GLT symbol 	given by $\mathbf{f}$. Every diagonal sampling sequence $\{D_\mathbf{n}(\mathbf{a})\}_n$ generated by a Riemann integrable function $\mathbf{a}(\boldsymbol{x})$, $\boldsymbol{x}\in  [0,1]^d$ is a GLT sequence and with GLT symbol given by $\mathbf{a}$.
	
  \item[GLT4] Every sequence which is distributed as the constant zero in the singular value sense is a GLT sequence with
	symbol~$0$.
	In particular:
\begin{itemize}
		\item every sequence in which the rank divided by the size tends to zero, as the matrix-size tends to infinity;
		\item every sequence in which the trace-norm (i.e., sum of the singular values) divided by the size tends to zero, as the matrix-size tends to infinity.
	\end{itemize}	
\end{description}

From a practical view-point, on the one hand, one of the main advantages for a sequence of belonging to the GLT class is that, under certain hypotheses, crucial spectral and singular value information can be derived using the concept of GLT symbol.
On the other hand, the above properties imply the following important features of the GLT symbol. Given a sequence $\{A_\mathbf{n}\}_\mathbf{n}$ obtained by algebraic operations  of a finite set of GLT sequences, the small-norm and low-rank  terms which composes the sequence should be neglected in the computation of the GLT symbol. Consequently, it happens that for small matrix-sizes $n$, the approximations may not be as accurate as it is desirable.

For this reason in \cite{Momentary_1} it has been introduced and exploited the concept of a (singular value and spectral) ``momentary symbols'', starting from a special case of Toeplitz structures. Here we generalize the notion to that of ``GLT momentary symbols'': the construction stems from that of the symbol in the GLT sense, but in practice the information of the small norm contributions is kept in the symbol and this may lead to higher accuracy, at least in some emblematic cases, when approximating the singular values and eigenvalues of Toeplitz-like matrices, even for small dimensions.

\subsection{The GLT momentary symbol sequence}
\label{sec:momentary}
For clarity in this subsection we consider the matrix-sequences in detail only in the unilevel and scalar setting. We want to avoid a cumbersome notation, but the ideas are extensible in a plain manner to the case where the involved GLT symbols are also matrix-valued and multivariate, as briefly sketched.

As an example, we take the following second-order differential equation with Dirichlet boundary conditions
\begin{equation*} \label{FD}
\begin{cases}
-(a(x)u'(x))' +b(x)u'(x) + c(x)u(x) = f(x), & x\in (0,1),\\
u(0) = \alpha, \qquad u(1) =\beta. &
\end{cases}
 \end{equation*}

The well-posedness of the previous diffusion-convection-advection problem holds in the case where $a(x)\in C^1(0,1)$. Furthermore, the uniqueness and existence of the  solution are guaranteed in the case where $a(x) > 0$, $c(x) \ge 0$ and with continuous functions $b(x),c(x)$ on $[0,1]$, with $f(x)\in L^2([0,1])$ (see \cite{Brezis}). For a more exhaustive discussion regarding the conditions of existence and uniqueness, even in the multidimensional case, we refer to \cite{Pozio,Punzo} and references therein.

 From a GLT viewpoint, we only require the following much weaker assumptions that are

\begin{itemize}
\item $a(x), c(x)$ are real-valued functions, continuous almost everywhere, defined in  $[0,1]$,
\item $b(x)$ is a real-valued function on $[0,1]$, such that $|b(x)x^{\alpha}|$ is bounded for some  $\alpha<3/2$,
\end{itemize}
while $f(x)$ is a general function.\\

We employ  central second-order finite differences  for approximating the given equation.
We  define the stepsize $h=\frac{1}{n+1}$ and the points $x_k=kh$  for $k$ belonging to the interval $[0,n+1]$. Let $a_k:= a(x_{\frac k2})$ for any $k\in[0,2n+2]$ and set $b_j:=b(x_j)$, $c_j:=c(x_j)$, $f_j:=f(x_j)$ for every $j=0,\ldots,n+1$. We compute approximations $u_j$ of the values $u(x_j)$ for $j=1,\ldots,n$ by solving the following linear system
\begin{equation}\label{linear-sys}
A_n\begin{pmatrix}
u_1\\
u_2\\
\vdots\\
u_{n-1}\\
u_n
\end{pmatrix} + B_n
\begin{pmatrix}
u_1\\
u_2\\
\vdots\\
u_{n-1}\\
u_n
\end{pmatrix}
+
C_n\begin{pmatrix}
u_1\\
u_2\\
\vdots\\
u_{n-1}\\
u_n
\end{pmatrix}
=h^2
\begin{pmatrix}
f_1 + \frac{1}{h^2}a_1\alpha + \frac{1}{2h}b_1\alpha\\
f_2\\
\vdots\\
f_{n-1}\\
f_n+ \frac{1}{h^2}a_{2n+1}\beta - \frac{1}{2h}b_{n}\beta
\end{pmatrix},
\end{equation}
where
\[
A_n =
\begin{pmatrix}
 a_1 + a_3  & - a_3 & & &  \\
 -a_3 & a_3+a_5 & -a_5 & &  \\
 & \ddots&\ddots &\ddots &  \\
 & &  -a_{2n-3}& a_{2n-3}+a_{2n-1}& -a_{2n-1} \\
 & & & -a_{2n-1}&a_{2n-1}+a_{2n+1}
\end{pmatrix},
\]
\[
B_n =
\frac{h}{2}
\begin{pmatrix}
 0  & b_1 & & & \\
 -b_2 & 0 & b_2 & & \\
 & \ddots&\ddots &\ddots &  \\
 & &  -b_{n-1} & 0 & b_{n-1} \\
 & & & -b_{n} & 0
\end{pmatrix},
\quad
C_n =
h^2
\diag(c_1,\ldots,c_n).
\]

In the case where $a(x)\equiv 1$ and $b(x)\equiv 1$, we find the basic Toeplitz structures
\begin{equation*}\label{basic toep1}
K_n = T_n(2-2\cos\theta)=
\begin{pmatrix}
 2  & - 1 & & &  \\
 -1 & 2 & -1 & &  \\
 & \ddots&\ddots &\ddots &  \\
 & &  -1& 2 & -1 \\
 & &     &  -1&  2
\end{pmatrix},
\end{equation*}
\begin{equation*}\label{basic toep2}
H_n = T_n(\mathbf{i}\sin\theta)=
\frac{1}{2}
\begin{pmatrix}
 0  & 1 & & & \\
 -1 & 0 & 1 & & \\
 & \ddots&\ddots &\ddots &  \\
 & &  -1 & 0 & 1 \\
 & & & -1 & 0
\end{pmatrix},
\end{equation*}
which are of importance since $A_n=D_n(a)K_n + E_n$ and $B_n=\textcolor{black}{h}D_n(b)H_n$ with $\{E_n\}_n\sim_\sigma 0$, as an immediate check in \cite{MR4157202} can show. Therefore, by using the GLT axioms (as done in detail in \cite{MR4157202}) we obtain
\[
\{ A_n\}_n \sim_{\textsc{glt}}
a(x)(2-2\cos\theta), \quad
\left\{ \frac{1}{h} B_n\right\}_n \sim_{\textsc{glt}}\mathbf{i} b(x)\sin\theta, \quad \left\{\frac{1}{h^2} C_n\right\}_n \sim_{\textsc{glt}} c(x).
\]
As a conclusion $\left\{B_n\right\}_n \sim_{\textsc{glt}}\text 0$, $\ \left\{C_n\right\}_n \sim_{\textsc{glt}} 0$ and hence, setting $X_n=A_n+B_n+C_n$ the actual coefficient matrix of the linear system in (\ref{linear-sys}), again by the $*$-algebra structure of the GLT matrix-sequences, we deduce
 \[
\{ X_n\}_n \sim_{\textsc{glt}}
a(x)(2-2\cos\theta).
\]
Now, following \cite{Momentary_1}, the idea is to consider not only the asymptotic setting, but also the case of moderate sizes. As a consequence, for increasing the precision of the evaluation of eigenvalues and singular values, we can associate to
\[
X_n=A_n+B_n+C_n
\]
the specific symbol $f_n(x,\theta)= a(x)(2-2\cos\theta)+h\mathbf{i} b(x)\sin\theta+ h^2c(x)$.

We are now in position to give a formal definition of GLT momentary symbols.

 \begin{defn}[GLT momentary symbols]
\label{def:momentarysymbols}
Let $\{X_n\}_n$ be a matrix-sequence and assume that there exist matrix-sequences $\{A_n^{(j)}\}_n$, scalar sequences $c_n^{(j)}$, $j=0,\ldots,t$,
and measurable functions $f_j$ defined over $[-\pi,\pi]\times [0,1]$, $t$ nonnegative integer independent  of $n$, such that
\begin{eqnarray} \nonumber
\left\{ \frac{A_n^{(j)}}{ c_n^{(j)}}\right\}_n &\sim_{\textsc{glt}} & f_j, \\ \nonumber
c_n^{(0)}=1, & & c_n^{(s)}=o(c_n^{(r)}), \ \ t\ge s>r, \\ \label{eq:sequence}
\{X_n\}_n & = & \{A_n^{(0)}\}_n + \sum_{j=1}^t \{A_n^{(j)}\}_n.
\end{eqnarray}
Then, by a slight abuse of notation,
\begin{equation}\label{eq:GLT_momentary_1D}
f_n=f_0+ \sum_{j=1}^t c_n^{(j)} f_j
\end{equation}

is defined as the GLT momentary symbol for $X_n$ and $\{f_n\}$ is the sequence of GLT momentary symbols for the matrix-sequence $\{X_n\}_n$.

\end{defn}

Of course, in line with Section \ref{sec:introduction:glt}, the momentary symbol could be matrix-valued with a number of variables equal to $2d$ and domain $[-\pi,\pi]^d\times [0,1]^d$ if the basic matrix-sequences appearing in Definition \ref{def:momentarysymbols} are, up to proper scaling, matrix-valued and multilevel GLT matrix-sequences.
\textcolor{black}{For example in the scalar $d$-variate setting relation (\ref{eq:GLT_momentary_1D}) takes the form
\[f_{\textbf{n}}= \sum_{\textbf{j}=\textbf{0}}^\textbf{t} c_\textbf{n}^{(\textbf{j})} f_\textbf{j},\]
which is a plain multivariate (possibly block) version of (\ref{eq:GLT_momentary_1D}).}

Clearly there is a link with the GLT theory stated in the next result.

\begin{thm}\label{moment-vs-glt}
Assume that the matrix-sequence $\{X_n\}_n$ satisfies the requirements in Definition \ref{def:momentarysymbols}. Then $\{X_n\}_n$ is a GLT matrix sequence and the GLT symbol $f_0$ of the main term $A_n^{(0)}$ is the GLT symbol of $\{X_n\}_n$, that is, $\{X_n\}_n \sim_{\textsc{glt}}  f_0$ and $
\lim_{n\to \infty} f_n=f_0$ uniformly on the definition domain.
\end{thm}

The given definition of momentary symbols is inspired, as it is clear from the initial example of diffusion-convection-advection equation, by the example of approximated differential equations, where the presence of differential operators of different orders induces, after a possible proper scaling, a structure like that reported in (\ref{eq:sequence}).

The idea is that the momentary symbol can be used for giving a more precise evaluation either of the spectrum or of the eigenvalues for moderate sizes of the matrices and not only asymptotically. However, we should be aware that, intrinsically, there is no general recipe especially for the eigenvalues. In fact, as already proven in \cite{MR2176808}, a rank one perturbation of infinitesimal spectral norm actually can change the spectra of matrix-sequences, sharing the same GLT symbol and even sharing the same sequence of momentary symbols.

\begin{description}
\item[Example 1] Take the matrices $T_n(e^{\mathbf{i} \theta})$ and $X_n=T_n(e^{\mathbf{i} \theta})+ e_1 e_n^T c_n^{(1)}$ with $c_n^{(1)}=n^{-\alpha}$, $\alpha>0$ any positive number independent of the matrix-size $n$. By direct inspection $\{e_1 e_n^T c_n^{(1)}\}_n\sim_\sigma 0$ and hence it is a GLT matrix-sequence with zero symbol, independently of the parameter $\alpha$. If we look at the GLT momentary symbols then they coincide with the GLT symbol for both $\{T_n(e^{\mathbf{i} \theta})\}_n$ and $\{X_n\}_n$: however while in the first case, the eigenvalues are all equal to zero, in the second case they distribute asymptotically as the GLT symbol $e^{\mathbf{i} \theta}$ (which is also the GLT momentary symbol for any $n$).
\item[Example 2]
Take a positive function $a$ defined on $[0,1]$ and the matrices $D_n(a)T_n(e^{\mathbf{i} \theta})$ and $X_n=D_n(a)T_n(e^{\mathbf{i} \theta})+ e_1 e_n^T c_n^{(1)}$ with $c_n^{(1)}=n^{-\alpha}$, $\alpha>0$ any positive number independent of the matrix-size $n$. Since $\{e_1 e_n^T c_n^{(1)}\}_n$ is a GLT matrix-sequence with zero symbol, independently of the parameter $\alpha$, we deduce that both $\{D_n(a)T_n(e^{\mathbf{i} \theta})\}_n$ and $\{X_n\}_n$ share the same GLT symbol $a(x)e^{\mathbf{i} \theta}$ (which is also the momentary symbol for any $n$). Again there is dramatic change: while in the first case, the eigenvalues are all equal to zero, in the second case they distribute asymptotically as the function $\hat a e^{\mathbf{i} \theta}$, where $\hat a$ is the limit (if it exists) of the geometric mean of sampling values present in $D_n(a)$, as $n$ tends to infinity: since $n^{-\alpha/n}$ converges to $1$ independently of the parameter $\alpha$ as $n$ tends to infinity, $\hat a$ will depend only on the diagonal values of $D_n(a)$. As a conclusion the eigenvalue distributions do not coincide with the GLT momentary symbols and this is a message that the present tool could be not effective and even misleading, when very non-normal matrices are considered.

In this setting it must be emphasized that the asymptotic eigenvalue distribution is discontinuous with respect to the standard norms or metrics widely considered in the context of matrix-sequences.
\end{description}
\section{\textcolor{black}{All-at-once solution of parabolic problems}}
\label{sec:problem}
The aim of this section is that of describing as accurate as possible the spectra and singular values of the structured linear system sequence stemming by the space-time discretization for a parabolic diffusion problem. Then, we consider the diffusion equation in one space dimension,
\begin{equation*}
u_t=u_{xx}, \quad x\in (a,b),\ t\in [0, T],
\end{equation*}
where we are  prescribing $u$ at $t= 0$ and imposing the periodicity condition $u(x\pm(b-a),t)=u(x,t)$.

We approximate our parabolic model problem on a rectangular space-time grid consisting of $N_t$ time intervals and $N_x$ space intervals. We obtain a sequence of linear systems, in which the each component is of the form
\begin{equation}
A_{\mathbf{n}}x=b, \quad A_\mathbf{n}=J_{N_t}\oplus Q_{N_x}=J_{N_t}\otimes  \mathbb{I}_{N_x}+\mathbb{I}_{N_t}\otimes Q_{N_x}\in\mathbb{R}^{N\times N}, \quad x,b\in\mathbb{R}^{N}, \label{eq:system}
\end{equation}
where $N=N_tN_x$, $\mathbf{n}=(N_t,N_x),$ $\mathbb{I}_{m}$ is the identity matrix of size $m$, and the matrices $J_{N_t}$ and $Q_{N_x}$ come from the discretization in time and space, respectively.
In the following, we describe the time and space discretization and, in particular, how this leads to structured components of the matrix $A_{\mathbf{n}}$.
\subsection{Time discretization}
\label{sec:problem:time}
The principal ingredients of the time discretization are:
\begin{itemize}
 \item Choosing $N_t$ equispaced points in $[0,T]$ with stepsize $h_t=T/N_t$, that is, $t_j=jh_t$, for $j=1,\ldots,N_t$.
 \item Discretizing in time by standard Euler backwards.
 \end{itemize}

Regarding notations, for the sake of simplicity, since we are considering a 2D problem, the symbols will have as Fourier variable $(\theta, \xi)$
instead of the standard choice $(\theta_1,\theta_2)$ indicated in the notations of Section \ref{sec:introduction:tep_circ} (see Definition \ref{def:toeplitz_block_generating}).

The resulting matrix is $J_{N_{t}}$, which has the following unilevel scalar Toeplitz structure:
  \begin{equation}
  J_{N_t}=\frac{1}{h_t}
\begin{bmatrix}
1 \\
-1 & 1\\
& \ddots & \ddots \\
& & -1 & 1
\end{bmatrix}=\frac{1}{h_t}T_{N_t}(f_J),\label{eq:problem:time:J}
  \end{equation}
 where $f_J$ is the generating function of the matrix-sequence $\{T_{n}(f_J)\}_{n}$ with
  \begin{equation*}
    f_J(\theta)=1-\E^{\mathbf{i}\theta}.\label{eq:problem:f_J}
  \end{equation*}

\subsection{Space discretization}
\label{sec:problem:space}
The principal elements of the time discretization are:
 \begin{itemize}
 \item Choosing $N_x$ equispaced points in $[a,b]$. Since we are considering periodic boundary conditions, we have step size $h_x=(a-b)/N_x$ and $x_j=h_x(j-1)$, for $j=1,\ldots,N_x$.
 \item Discretizing in space using second order finite differences.
 \end{itemize}
  Consequently, the space discretization matrix will be the circulant matrix $Q_{N_x}$ of the form:
    \begin{equation*}
    Q_{N_x}=
    \frac{1}{h_x^2}
    \begin{bmatrix}
    2 &-1&&&-1\\
    -1 & 2& -1\\
    & \ddots & \ddots &\ddots \\
    & & -1 & 2 & -1\\
    -1& & & -1 & 2
    \end{bmatrix}= \frac{1}{h_x^2}C_{N_x}(f_{Q}),\label{eq:problem:space:Q}
    \end{equation*}
    where
    \begin{equation*}
    f_Q(\xi)=2-2\cos\xi\label{eq:problem:f_Q}
    \end{equation*}
 is the   generating  function of the matrix.

Of course a different choice of Dirichlet boundary conditions would lead to the standard discrete Laplacian $T_{N_x}(f_Q)$: the analysis is equivalent since also this matrix admits a well known diagonalization matrix, that is the sine transform matrix of type I, which is real, orthogonal and symmetric.

\subsection{Analysis of the coefficient matrix $A_\mathbf{n}$}
\label{sec:coefficientmatrix}
We have seen that discretizing of the problem of interest for a sequence of discretization parameters $h_x$ and $h_t$ leads to a sequence of linear systems, whose approximation error tends to zero as the coefficient matrix-size grows to infinity. The $\mathbf{n}$th coefficient matrix component is of the form
\begin{equation}\label{to be diagonalized}
A_\mathbf{n}=\frac{1}{h_t}T_{N_t}(f_J)\otimes \mathbb{I}_{N_x}+\mathbb{I}_{N_t}\otimes \frac{1}{h_x^2}C_{N_x}(f_Q).
\end{equation}

In order to design efficient solvers for the considered linear systems, it is of crucial importance to know the spectral propriety of the matrix-sequence $\{A_\mathbf{n}\}_\mathbf{n}$.
Hence, this section is devoted to the analysis of the structure of the matrix-sequence  $\{{A}_\mathbf{n}\}_\mathbf{n}$  in  \eqref{eq:system}.
In particular, we  provide the singular values and spectral analysis using algebraic tricks, the GLT theory, and the concept of GLT momentary symbols.
\subsection{GLT analysis of the coefficient sequence $\{{A}_\mathbf{n}\}_{\mathbf{n}}$}\label{sec:general_glt_analysis}
The asymptotic spectral and singular value distribution, for the matrix-size $d(\mathbf{n})$ sufficiently large, of the matrix-sequence $\{{A}_\mathbf{n}\}_\mathbf{n}$
depend on how $h_x$ and $h_t$ approaches zero. Let $c_h\coloneqq h_x^2/h_t$, we have three different cases to consider.
\begin{itemize}
\item[\textsc{Case} 1.] $\left[c_h\to\infty\right]:$ If $h_t\to 0$ faster than $C_1h_x^2$, where $C_1$ is a constant, then we can consider the matrix
\begin{equation*}
h_t{A}_\mathbf{n}=T_{N_t}(f_J)\otimes \mathbb{I}_{N_x}+\mathbb{I}_{N_t}\otimes \underbrace{\frac{h_t}{h_x^2}}_{c_h^{-1}\to 0}C_{N_x}(f_Q).
\end{equation*}
Then, the sequence $\{h_t A_\mathbf{n}\}_\mathbf{n}=\{T_{N_t}(f_J)\otimes \mathbb{I}_{N_x}+\mathbb{N}_\mathbf{n}\}_\mathbf{n}$, where $\mathbb{N}_\mathbf{n}$ is a small-norm matrix in the sense of the item 2 of property \textbf{GLT4}, with $\|\mathbb{N}_\mathbf{n}\|<C_2$, $C_2$ constant. Consequently, from \textbf{GLT4}, $\{\mathbb{N}_\mathbf{n}\}_\mathbf{n}$ is a matrix-sequence distributed in the singular value sense as $0$, which implies that $\{\mathbb{N}_\mathbf{n}\}_\mathbf{n}$ is zero-distribued in GLT sense as described in \textbf{GLT4}. Moreover $f_J$ is a trigonometric polynomial, then Theorem \ref{thm:prel_toepl_szeg_multi}, properties \textbf{GLT1}-\textbf{GLT4} and Lemma \ref{lemm:tensor_prod} imply that
\begin{equation*}
\{h_t A_\mathbf{n}\}_\mathbf{n}\sim_{\textsc{glt}}f_J(\theta)\otimes 1+1\otimes0= f_A^{(1)}(\theta,\xi).
\end{equation*}

The $1$ present in $ f_J(\theta)\otimes 1$ should be interpreted as $1\E^{0\mathbf{i}\xi}$, and $1\otimes0$ should be interpreted as $1\E^{0\mathbf{i}\theta}\otimes0\E^{0\mathbf{i}\xi}$. Hence, the GLT symbol of the sequence $ \{h_t A_\mathbf{n}\}_\mathbf{n}$
is the bivariate function
\begin{equation*}
f_A^{(1)}(\theta,\xi)= f_J(\theta)=1-\E^{\mathbf{i}\theta},\label{eq:coefficientmatrix:fA1}
\end{equation*}
and it should be interpreted as the function $f_J(\theta)\otimes 1$, with is constant in the second component.

From the property \textbf{GLT1}, the function $f_A^{(1)}(\theta,\xi)$ describes the singular value distribution in the sense of relation (\ref{eq:introduction:background:distribution:sv}). More in detail
\begin{equation*}
\{h_t A_\mathbf{n}\}_\mathbf{n}\sim_{\textsc{glt},\sigma}1-\E^{\mathbf{i}\theta}.
\end{equation*}

 However, the matrix-sequence $\{h_t A_\mathbf{n}\}_\mathbf{n}$ is not symmetric, hence the distribution does not hold in the eigenvalue sense
 (see also {\bf Example 1} and {\bf Example 2} at the end of Section \ref{sec:momentary}). Because of the structure of $J_{N_t}=\frac{1}{h_t}T_{N_t}(f_J)$ in equation (\ref{eq:problem:time:J}) it is straightforward to see that the asymptotic spectral distribution is given by ${\mathfrak{f}({\theta},\xi)}=1,$ accordingly to relation (\ref{eq:introduction:background:distribution:ev}), that is
\begin{equation*}
\{h_t A_\mathbf{n}\}_\mathbf{n}\sim_{\lambda}1.
\end{equation*}
\item[\textsc{Case} 2.]  $\left[c_h\to0\right]:$ If  $h_x^2\to 0$ faster than $C_1h_t$, where $C_1$ is a constant, then we have
\begin{equation*}
h_x^2A_\mathbf{n}=\underbrace{\frac{h_x^2}{h_t}}_{c_h\to 0}T_{N_t}(f_J)\otimes \mathbb{I}_{N_x}+\mathbb{I}_{N_t}\otimes C_{N_x}(f_Q).
\end{equation*}
Then, the sequence $\{h_x^2 A_\mathbf{n}\}_\mathbf{n}=\{\mathbb{N}_\mathbf{n}+\mathbb{I}_{N_t}\otimes C_{N_x}(f_Q)\}_\mathbf{n}$, where $\mathbb{N}_\mathbf{n}$ is a small-norm matrix in the sense of the item 2 of property \textbf{GLT4}, with $\|\mathbb{N}_\mathbf{n}\|<C_2$, $C_2$ constant.
Then, $\{\mathbb{N}_\mathbf{n}\}_\mathbf{n}$ is a matrix-sequence distributed in the singular value sense, and consequently in the GLT sense, as $0$. Moreover, ${f_Q}$ belongs to the Dini-Lipschitz class, consequently, properties \textbf{GLT2}-\textbf{GLT4},  and Lemma \ref{lemm:tensor_prod}  imply that
\begin{equation*}
\{h_x^2 A_\mathbf{n}\}_\mathbf{n}\sim_{\textsc{glt}}0\otimes 1+1\otimes f_Q(\xi)=1\otimes f_Q(\xi)=f_A^{(2)}(\theta,\xi),
\end{equation*}
where the GLT symbol is given by
\begin{equation*}
 f_A^{(2)}(\theta,\xi)= f_Q(\xi)=2-2\cos\xi.\label{eq:coefficientmatrix:fA2}
\end{equation*}
In this case the function $f_A^{(2)}(\theta,\xi)$ is a singular value symbol for the sequence $\{h_x^2 A_\mathbf{n}\}_\mathbf{n}$, and also an eigenvalue symbol, since the matrices $C_{N_x}(f_{Q})$ are Hermitian for each $N_x$. Hence we have
\begin{equation*}
\{h_x^2 A_\mathbf{n}\}_\mathbf{n}\sim_{\textsc{glt},\sigma, \lambda} 2-2\cos\xi.
\end{equation*}
\item[\textsc{Case} 3.] $\left[c_h=c= \text{constant}\right]:$ The last case is when $h_x^2$ and $h_t$ are proportional  and related by the constant $c_h=c=\frac{h_x^2}{h_t}$, independent of the various step-sizes.
In this setting we have
\begin{equation*}
h_x^2A_\mathbf{n}=\underbrace{\frac{h_x^2}{h_t}}_{c_h}T_{N_t}(f_J)\otimes \mathbb{I}_{N_x}+\mathbb{I}_{N_t}\otimes C_{N_x}(f_Q).
\end{equation*}
Consequently, from \textbf{GLT2}, \textbf{GLT3} and Lemma \ref{lemm:tensor_prod},  the following relationship holds when $c_h$ is a constant,
\begin{equation*}
\{h_x^2A_\mathbf{n}\}_\mathbf{n}\sim_{\textsc{glt}}c f_J(\theta)\otimes 1+1\otimes f_Q(\xi)=f_A^{(3)}(\theta,\xi).
\end{equation*}
From considerations analogous to the case 1 and 2 we have
\begin{equation*}
\{h_x^2A_\mathbf{n}\}_\mathbf{n}\sim_{\textsc{glt},\sigma}c(1-\E^{\mathbf{i}\theta})+(2-2\cos\xi).\label{eq:coefficientmatrix:fA3}
\end{equation*}
Since the matrix $h_x^2A_\mathbf{n}$ is not Hermitian, the eigenvalue symbol $\es({\theta},\xi)$ cannot be directly derived by $f_A^{(3)}(\theta,\xi)$
(see again the discussion in the examples after Definition \ref{def:momentarysymbols}).

In this setting the situation is simple because the involved twolevel structure can be simply block-diagonalized, while the use of the GLT momentary symbol becomes useful in approximation the singular values of the sequence $\{h_x^2A_\mathbf{n}\}_\mathbf{n}$.
\end{itemize}

\subsection{Analysis of the coefficient matrix-sequence  $\{{A}_\mathbf{n}\}_{\mathbf{n}}$ by algebraic manipulations and GLT momentary symbols}\label{sec:general_momentary_analysis}

The first observation is that the matrix in (\ref{to be diagonalized}) admits a perfect decomposition which shows in evidence a lower triangular matrix, which is similar to the original one and hence all the eigenvalues are known exactly.
In fact, by looking carefully at (\ref{to be diagonalized}), we obtain that
\[
\frac{1}{h_t}T_{N_t}(f_J)\otimes \mathbb{I}_{N_x}= \mathbb{I}_{N_t}\left[\frac{1}{h_t}T_{N_t}(f_J)\right]\mathbb{I}_{N_t}\otimes \mathbb{F}_{N_x}\mathbb{I}_{N_x}
\mathbb{F^*}_{N_x}
\]
and
\[
\mathbb{I}_{N_t} \otimes \frac{1}{h_x^2}C_{N_x}(f_Q)= \mathbb{I}_{N_t}\mathbb{I}_{N_t}\mathbb{I}_{N_t} \otimes \mathbb{F}_{N_x}\frac{1}{h_x^2}D_{N_x}\mathbb{F^*}_{N_x},
\]
where $\mathbb{F}_{N_x}$ is the unitary Fourier matrix of size $N_x$, $\mathbb{F^*}_{N_x}$ is its transpose conjugate and hence its inverse, and
$D_{N_x}$ is the diagonal matrix containing the eigenvalues of $C_{N_x}(f_Q)$ that is $f_Q(2\pi j/N_x)=2-2\cos(2\pi j/N_x)$, $j=0,1,\ldots, N_x-1$.

Since $T_{N_t}(f_J)$ is lower bidiagonal matrix with $1$ on the main diagonal, it can be easily seen that the eigenvalues of $A_\mathbf{n}$ in (\ref{to be diagonalized}) are exactly
\[
\frac{1}{h_t}+\frac{1}{h_x^2}(2-2\cos(2\pi j/N_x)),\ \ \ j=0,1,\ldots, N_x-1,
\]
each of them with multiplicity $N_t$. As a consequence, by taking a proper normalization, the spectral radius $\rho(h_x^2 A_{\mathbf{n}})$ will coincide simply with $4+c_h$.

It is clear that, in this context, due to the high non-normality of the term $T_{N_t}(f_J)$, after proper scalings depending on $h_t$ and $h_x$, the eigenvalues are a uniform sampling of a function which is not the GLT symbol and is not the associated GLT momentary symbol.
This is not surprising given the discussion regarding the asymptotical behaviour of the matrix-sequences reported in {\bf Example 1} and in {\bf Example 2}, when discussing the potential and the limitations of the notion of GLT momentary symbols.

Also in this setting, by imposing (quite artificial) periodic boundary conditions in time, the term $T_{N_t}(f_J)$ will change into $C_{N_t}(f_J)$ and magically a one-rank correction repeated $N_t$ times to the matrix $A_\mathbf{n}$ will produce a new matrix with the same GLT and momentary symbols as before: however in this case the eigenvalues will be exactly the sampling of such functions. This is a further confirmation of the delicacies of the eigenvalues that can have dramatic changes due to minimal corrections, when we are in a context of higlhy non-normal matrices.

\subsubsection{Singular values of $h_x^2A_\mathbf{n}$ (exact)}

The singular values $\sigma_{1}(h_x^2A_\mathbf{n}),\dots,$ $\sigma_{d(\mathbf{n})}(h_x^2A_\mathbf{n}) $ of the matrix $h_x^2A_\mathbf{n}$ are given the square root of the eigenvalues of the Hermitian matrix $h_x^4A_\mathbf{n}A_\mathbf{n}^{\textsc{t}}$. Hence, in order to provide exactly $\sigma_{i}(h_x^2A_\mathbf{n})$, $i=1\dots, d(\mathbf{n})$, we are interested at the spectrum of the matrix
\begin{equation*}
\begin{split}
&h_x^4A_\mathbf{n}A_\mathbf{n}^{\textsc{t}}=\\
&\left[
\begin{array}{cccccccccc}
\tilde{Q}_{N_x}^2&-c_h\tilde{Q}_{N_x}\\
-c_h\tilde{Q}_{N_x}&\tilde{Q}_{N_x}^2+c_h^2\mathbb{I}_{N_x}&-c_h\tilde{Q}_{N_x}\\
   &-c_h\tilde{Q}_{N_x}&\tilde{Q}_{N_x}^2+c_h^2\mathbb{I}_{N_x}&\ddots&\\
   & & & &\\
   &   &\ddots&\ddots&-c_h\tilde{Q}_{N_x}\\
   &      & &-c_h\tilde{Q}_{N_x}&\tilde{Q}_{N_x}^2+c_h^2\mathbb{I}_{N_x}
\end{array}
\right],
\end{split}
\end{equation*}
where $\tilde{Q}_{N_x}=C_{N_x}+c_h\mathbb{I}_{N_x} $. Note that  $h_x^4A_\mathbf{n}A_\mathbf{n}^{\textsc{t}}$ is not a pure block-tridiagonal Toeplitz because of the missing constant $c_h^2$ in the block in the top left corner.
However,  for each fixed $N_t$ and $N_x$, the matrix $\tilde{Q}_{N_x}$ is a circulant matrix with generating function $f_{\tilde{Q}_{N_x}}(\xi)=2-2\cos\xi +c_h$, which is also its GLT momentary symbol. Thus we infer that $h_x^4A_\mathbf{n}A_\mathbf{n}^{\textsc{t}}$ is similar to a matrix $X_{d(\mathbf{n})}$, whose explicit expression is reported below
\begin{equation*}
\begin{split}
&
h_x^4A_\mathbf{n}A_\mathbf{n}^{\textsc{t}}\sim X_{d(\mathbf{n})}=\\
&\left[
\begin{array}{cccccccccc}
D_{\tilde{Q}}^{2}&-c_hD_{\tilde{Q}}\\
-c_hD_{\tilde{Q}}&D_{\tilde{Q}}^{2}+c_h^2\mathbb{I}_{N_x}&-c_hD_{\tilde{Q}}\\
   &-c_hD_{\tilde{Q}}&D_{\tilde{Q}}^{2}+c_h^2\mathbb{I}_{N_x}&\ddots\\
   &   &\ddots&\ddots&-c_hD_{\tilde{Q}}\\
   &   &      &-c_hD_{\tilde{Q}}&D_{\tilde{Q}}^{2}+c_h^2\mathbb{I}_{N_x}
\end{array}
\right],
\end{split}
\end{equation*}
with $D_{\tilde{Q}}=\diag_{\ell=1,\dots, N_x}\left(f_{\tilde{Q}_{N_x}}(\xi_{\ell,N_x}) \right)$.  Consequently we study the spectrum of $X_{d(\mathbf{n})}$ to attain formulas for the exact  singular values of $h_x^2A_{\mathbf{n}}$. Let us consider a permutation matrix $P$ such transforms $X_{d(\mathbf{n})}$ into an $N_t\times N_t$ block diagonal matrix $PX_{d(\mathbf{n})}P^{\textsc{t}}$, which has on the main diagonal, for $k=1,\ldots,N_x$,   blocks of the form

\begin{equation}
\begin{split}
&\left(PX_{d(\mathbf{n})}P^{\textsc{t}}\right)_{i,j=(k-1)N_t+1}^{\textcolor{black}{kN_t}}=\\
&\left[
\begin{array}{cccccccccc}
C_k^2&-c_hC_k\\
-c_hC_k&C_k^2+c_h^2&-c_hC_k\\
   &-c_hC_k&C_k^2+c_h^2&\ddots\\
   &   &\ddots&\ddots&-c_hC_k\\
   &   &      &-c_hC_k&C_k^2+c_h^2
\end{array}
\right],\label{eq:coefficientmatrix:submatrix}
\end{split}
\end{equation}
where $C_k=D_{\tilde{Q}}(k,k)=f_{\tilde{Q}_{N_x}}(\xi_{k,N_x})$.

Hence, the  union of the eigenvalues of all blocks of $PX_{d(\mathbf{n})}P^{\textsc{t}}$ is equivalent to the full spectrum of $h_x^4A_\mathbf{n}A_\mathbf{n}^{\textsc{t}}$. These local eigenvalue problems can be solved analytically (or numerically) independently from each other. For example for $N_t=2$ we have for every $k=1,\ldots,N_x$ the characteristic equation
\begin{equation*}
\left|
\begin{array}{cccccccccc}
\left(f_{\tilde{Q}_{N_x}}(\xi_{k,N_x})\right)^2-\lambda&-c_h(f_{\tilde{Q}_{N_x}}(\xi_{k,N_x}))\\
-c_h(f_{\tilde{Q}_{N_x}}(\xi_{k,N_x}))&\left(f_{\tilde{Q}_{N_x}}(\xi_{k,N_x})\right)^2+c_h^2-\lambda
\end{array}
\right|=0.
\end{equation*}
Thus, we have as singular value the union for $k=1,\dots, N_x$ of the quantities
\begin{equation*}
\begin{split}
\sigma^{(1)}(k,c_h)&=\sqrt{\frac{2(f_{\tilde{Q}_{N_x}}(\xi_{k,N_x}))^2+c_h^2}{2}-\frac{c_h}{2}\sqrt{4(f_{\tilde{Q}_{N_x}}(\xi_{k,N_x}))^2+c_h^2}},\\
\sigma^{(1)}(k,c_h)&=\sqrt{\frac{2(f_{\tilde{Q}_{N_x}}(\xi_{k,N_x}))^2+c_h^2}{2}+\frac{c_h}{2}\sqrt{4(f_{\tilde{Q}_{N_x}}(\xi_{k,N_x}))^2+c_h^2}}.
\end{split}
\end{equation*}
Clearly, solving the characteristic equation for $k=1,\ldots,N_x$ becomes more and more complex as $N_t$ grows. Hence, in the next section provide two possible approximations given by the GLT theory and by the GLT momentary formulations.
\subsubsection{Singular values of $h_x^2A_{\mathbf{n}}$ (approximation) via GLT momentary symbols}
For case (2),  in Section \ref{sec:general_glt_analysis} and in Section \ref{sec:general_momentary_analysis}, we have already shown that
\begin{equation*}
\{h_x^2 A_\mathbf{n}\}_\mathbf{n}\sim_{\sigma}\s_{A}^{(2)}(\theta,\xi)= 2-2\cos\xi.
\end{equation*}
\textcolor{black}{On the other hand, 
the subsequent sequence $\{f_\mathbf{n}^{(2)}\}_{\textbf{n}}$ with
\begin{equation*}
f_\mathbf{n}^{(2)}(\theta,\xi)= c_h(1-\E^{\mathbf{i}\theta})+(2-2\cos\xi),
\end{equation*}
is the sequence of GLT momentary functions.
}

 Remark \ref{rem:introduction:background:1} suggests to exploit these relations in order to obtain a better approximation of the singular value of $h_x^2A_{\mathbf{n}}$, with respect to the information obtained by the pure GLT symbol. In the following, we compute the quantities $|\s_A^{(2)}(\theta,\xi)|$ and $|f_\mathbf{n}^{(2)}(\theta,\xi)|$ using the specific grid described below
\begin{equation}
\theta_{j,N_t}=\frac{j\pi}{N_t\textcolor{black}{+1}}, \, j=1,\dots, N_t\qquad \xi_{\ell,N_x}=\frac{2\pi(\ell-1)}{N_x}\, \ell=1,\dots, N_x.
\label{eq:ANperfectgrids}
 \end{equation}
 \textcolor{black}{
We can observe in Figure \ref{fig:sigma} that the singular value of $h_x^2A_{\mathbf{n}}$ (blue circles) are well approximated by the
samplings of $|f_\mathbf{n}^{(2)}(\theta,\xi)|$ on the grid \eqref{eq:ANperfectgrids} (red stars). The
approximation by using  $|\s_A^{(2)}(\theta,\xi)|$, instead, is good when $c_h$ is small, see the top panel of Figure
\ref{fig:sigma} for $N_t=2$ and $N_x=10$, but it tends to become a substantially less accurate approximation otherwise, see the bottom panel of Figure \ref{fig:sigma} and Figure \ref{fig:error_sigma} where $N_t=10$ and $N_x=10$.}

\begin{figure}[!h]
\centering
\includegraphics[width=0.7\textwidth]{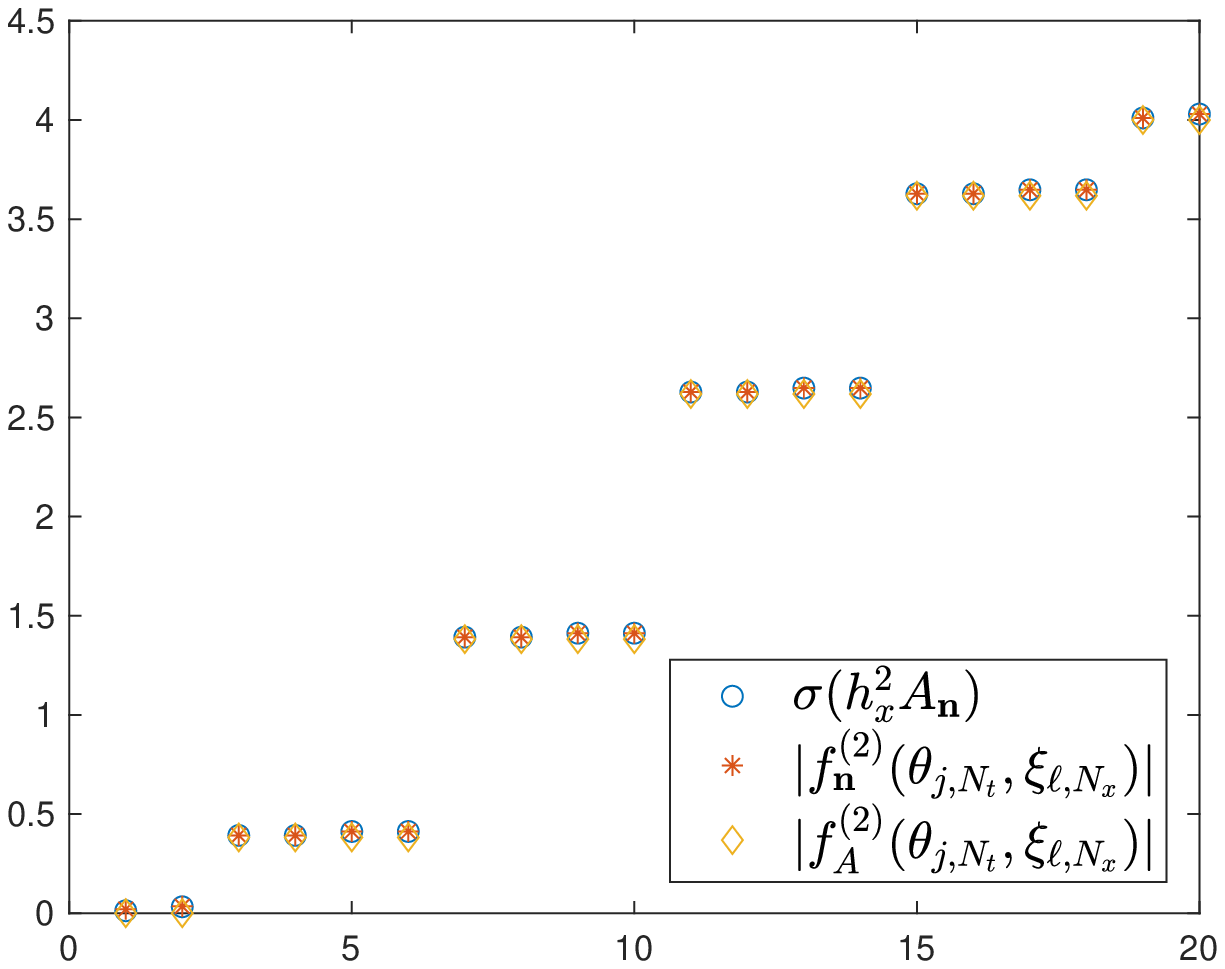}
\includegraphics[width=0.7\textwidth]{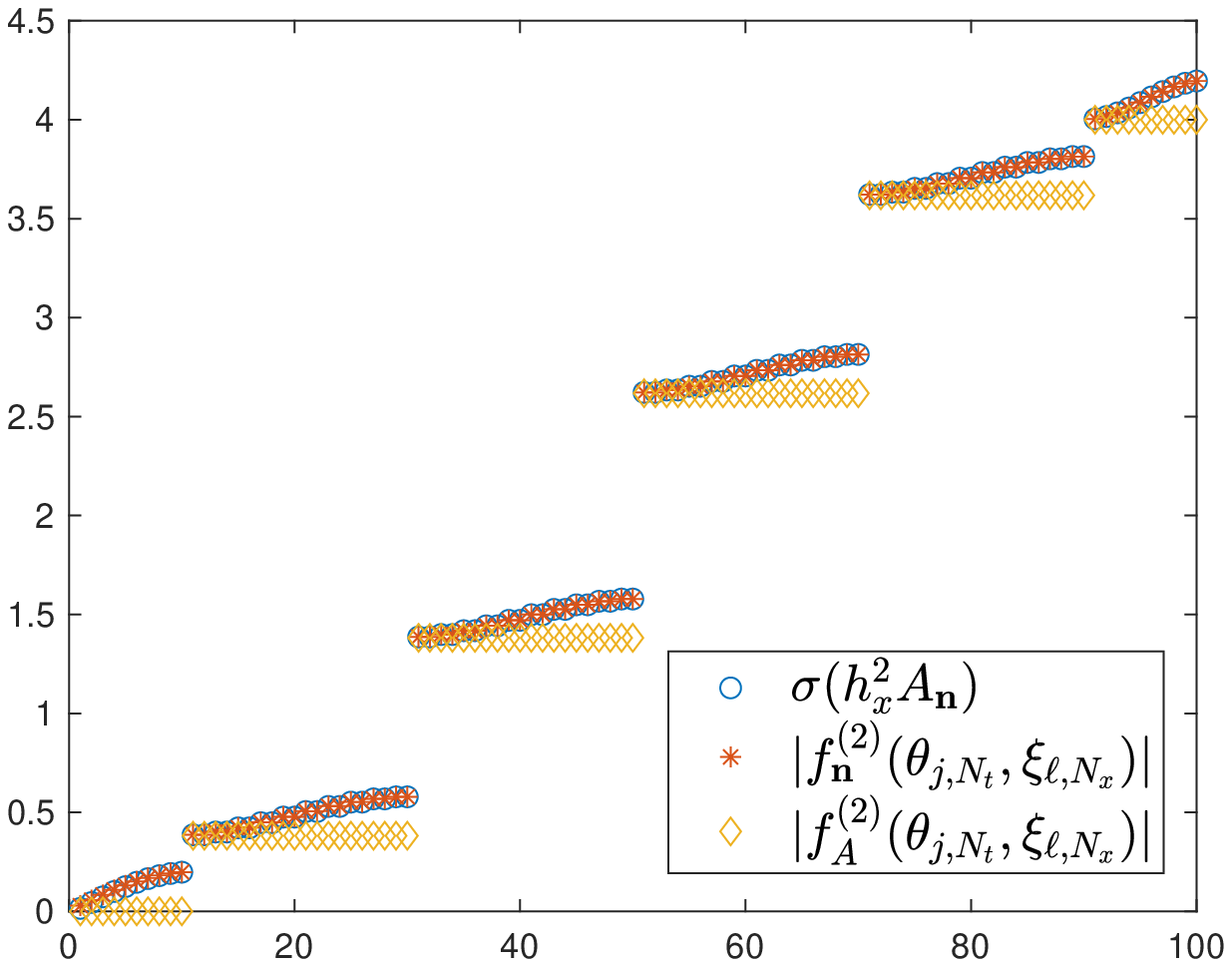}
\caption{\textcolor{black}{Singular values of $h_x^2A_{\mathbf{n}}$  and samplings of $|f_A^{(2)}(\theta,\xi)|$ and $|f_\mathbf{n}^{(2)}(\theta,\xi)|$ on the grid  \eqref{eq:ANperfectgrids}, for $N_t=2$ and $N_x=10$ (top) and $N_t=10$ and $N_x=10$ (bottom). }}
\label{fig:sigma}
\end{figure}

\begin{figure}[!h]
\centering
\includegraphics[width=1\textwidth]{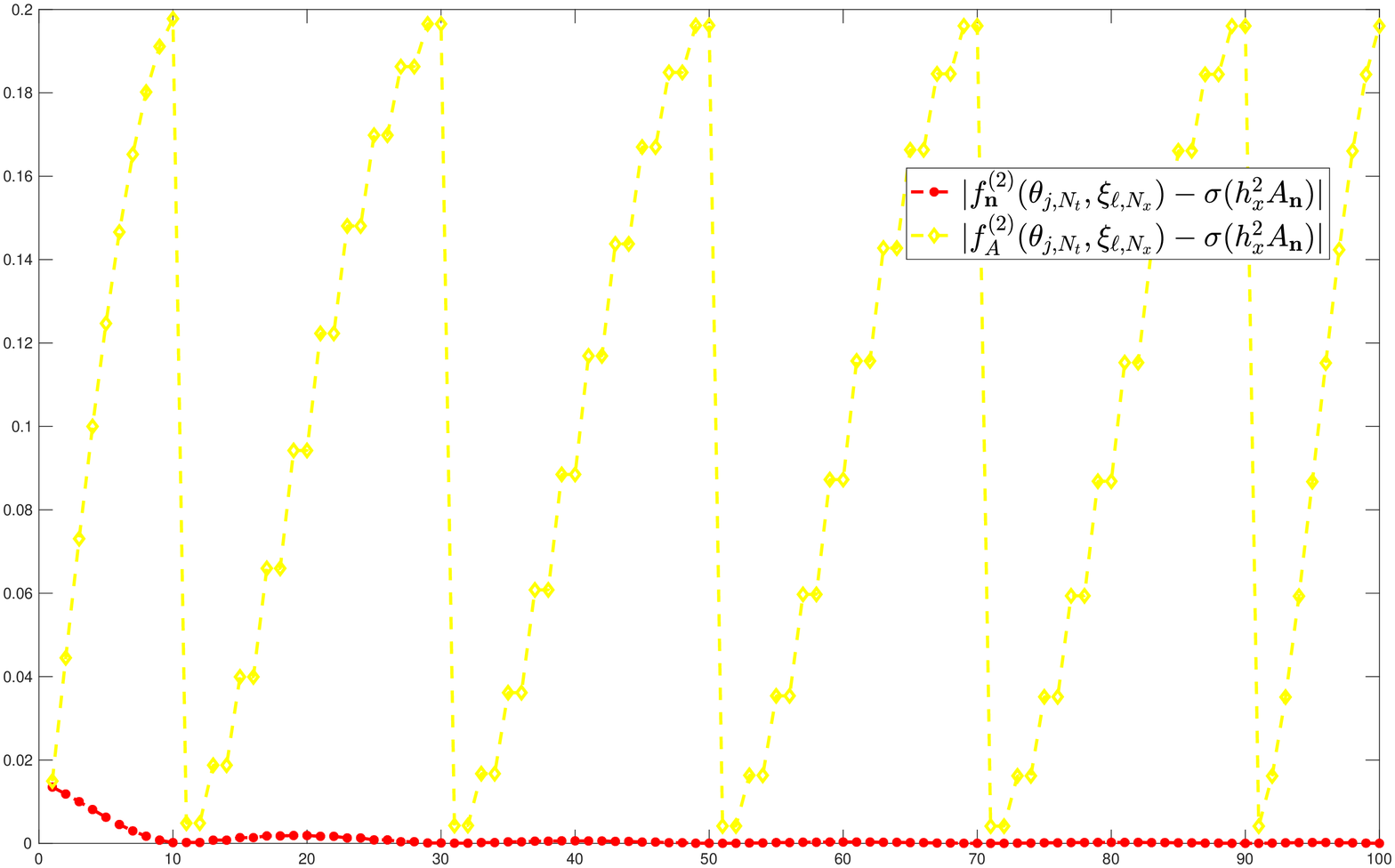}
\caption{\textcolor{black}{Singular values and samplings of $|f_A^{(2)}(\theta,\xi)|$ and $|f_\mathbf{n}^{(2)}(\theta,\xi)|$ for $N_t=10$ and $N_x=10$ on the grid in \eqref{eq:ANperfectgrids}.}}\label{fig:error_sigma}
\end{figure}
\subsubsection{2-norm of $h_x^2A_{\mathbf{n}}$ (approximation)}
In the following we are interested in providing a bound for the $2$-norm of the matrix $h_x^2A_{\mathbf{n}}$. By definition it is given by $\|h_x^2A_{\mathbf{n}}\|_2=\max_{j=1\dots,d(\mathbf{n})}|\sigma_{j}(h_x^2A_{\mathbf{n}})|$.
From the previous section we know that it can be computed by making the square root of the maximum eigenvalue of the block  in \eqref{eq:coefficientmatrix:submatrix}, corresponding to $\xi_{N_x/2+1,N_x}=\pi$. Since the $ \max_k \left( f_{\tilde{Q}_{N_x}}(\xi_{k,N_x})\right)=4+c_h$, we are interested in estimate the maximum eigenvalue of
\begin{equation}\label{eq:symmetriz_4_sing}
\left[
\begin{array}{cccccccccc}
(4+c_h)^2&-c_h(4+c_h)\\
-c_h(4+c_h)&(4+c_h)^2+c_h^2&-c_h(4+c_h)\\
   &-c_h(4+c_h)&(4+c_h)^2+c_h^2&\ddots\\
   &   &\ddots&\ddots&-c_h(4+c_h)\\
   &   &      &-c_h(4+c_h)&(4+c_h)^2+c_h^2
\end{array}
\right].
\end{equation}
 For this purpose we exploit the  concept of {$\tau_{\varepsilon,\varphi}$}-algebras of Subsection \ref{sec:matrix_algebra}.

In our case $$a=(4+c_h)^2+c_h^2, \quad b=-c_h(4+c_h),$$ and the matrix belongs to the $\tau_{\frac{c_h}{4+c_h},0}$-algebra, since the element with indices $i,j=1$ is $a+(c_h/(4+c_h))b$.
% \begin{align}
% \underbrace{
% \frac{(n-1)\pi}{n}}_{\tau_{11}}
% <\underbrace{
% \frac{(n-1/2)\pi}{n}}_{\tau_{\pm 1\mp1}}
% <\underbrace{
% \frac{n\pi}{n}}_{\tau_{-1-1}}
% <\underbrace{
% \frac{(n+1/2)\pi}{n}}_{\tau_{\pm 1\mp1}}
% \end{align}

Hence, we have $g(\theta)=(4+c_h)^2+c_h^2-2c_h(4+c_h)\cos\theta$, which coincides with the  eigenvalue symbol $\mathfrak{g}_{\mathbf{n}}$ of the matrix (\ref{eq:symmetriz_4_sing}).  Due to the Interlacing theorem, see~\cite{interlacing} and the specific relation between algebras \cite{Momentary_1}, the following relationships can be derived
\begin{align}
&\underbrace{\mathfrak{g}_{\mathbf{n}}\left(\frac{\pi(N_t-1/2)}{N_t+1/2}\right)}_{\mathrm{max}(\lambda_j(T_{N_t,1,0}(g)))}<
\underbrace{\|h_x^2A_{\mathbf{n}}\|_2^2}_{\mathrm{max}(\lambda_j(T_{N_t,c_h/(4+c_h),0}(g)))}%\mathrm{max}(\tau(c_h/(4+c_h),0))}
<\\
&
\underbrace{\mathfrak{g}_{\mathbf{n}}\left(\frac{\pi N_t}{N_t+1}\right)}_{\mathrm{max}(\lambda_j(T_{N_t,0,0}(g)))}
% &<
% \underbrace{\mathfrak{g}_{\mathbf{n}}\left(\frac{\pi N_t}{N_t+1/2},c_h\right)}_{\mathrm{max}(\tau(-1,0))}\\
<
\underbrace{\mathfrak{g}_{\mathbf{n}}(\pi)}_{\mathrm{max}(\lambda_j(T_{N_t,-1,-1}(f)))=\mathrm{max}(\mathfrak{g}_{\mathbf{n}})}.\nonumber
\end{align}
As a consequence, good upper and lower bounds for the 2-norm of $h_x^2A_{\mathbf{n}}$ are reported in the following set of inequalities
\begin{equation*}
\sqrt{\mathfrak{g}_{\mathbf{n}}\left(\frac{\pi(N_t-1/2)}{N_t+1/2}\right)}
<
\|h_x^2A_{\mathbf{n}}\|_2
<
\sqrt{\mathfrak{g}_{\mathbf{n}}\left(\frac{\pi N_t}{N_t+1}\right)}.
\end{equation*}
In Table~\ref{tbl:2normAN} we present approximations of the 2-norm of $h_x^2A_{\mathbf{n}}$ using the grid sampling from $\tau_{1,0}$ (lower bound), $\tau_{0,0}$ (upper bound), and $\tau_{-1,-1}$. Note that the sampling on the latter grid  is equivalent  to do the sampling of the singular value momentary symbols $f_\mathbf{n}^{(2)}(\theta,\xi)$ at their maximum point. The two-norm $\|h_x^2A_\mathbf{n}\|_2$ is computed numerically. We see that the 2-norm is well described by the two bounds given above, as $N_t$ increases. Hence, for this type of examples, the GLT momentary symbols provide, at least for moderate sizes, a more precise alternative to the pure GLT symbol.
% \begin{figure}[!ht]
% \centering
% \includegraphics[width=\textwidth]{exmp1.png}
% \caption{$N_t=10$, $c_h=1$, $g(\theta;c_h)=(4+c_h)^2+c_h^2-2(c_h(4+c_h))\cos\theta$, $g(\theta;1)=26-10\cos\theta$}
% \label{fig:grids}
% \end{figure}
\begin{table}[!h]
\centering
\caption{Approximations of $2$-norm for different $N_t$ and $c_h$. The maximum bound is $\sqrt{\max \mathfrak{g}_{\mathbf{n}}}=4+2c_h$.}
\label{tbl:2normAN}
\begin{tabular}{rrrrrr}
\toprule
$N_t$&$c_h$&$\sqrt{\mathfrak{g}_{\mathbf{n}}\left(\frac{\pi(N_t-1/2)}{N_t+1/2}\right)}$&$\|h_x^2A_\mathbf{n}\|_2$&$\sqrt{\mathfrak{g}_{\mathbf{n}}\left(\frac{\pi N_t}{N_t+1}\right)}$&$4+2c_h$\\
\midrule
1&1/8&4.06394205&4.12500000&4.12689350&4.25\\
10&1/8&4.24460651&4.24505679&4.24508270&4.25\\
100&1/8&4.24994073&4.24994128&4.24994131&4.25\\
1000&1/8&4.24999940&4.24999940&4.24999940&4.25\\
\midrule
1&1&4.58257569&5.00000000&5.09901951&6.00\\
10&1&5.96286240&5.96511172&5.96614865&6.00\\
100&1&5.99959287&5.99959555&5.99959689&6.00\\
1000&1&5.99999589&5.99999589&5.99999590&6.00\\
\midrule
1&8&10.58300524&12.00000000&14.42220510&20.00\\
10&8&19.78560029&19.78964627&19.80461186&20.00\\
100&8&19.99765486&19.99765952&19.99767802&20.00\\
1000&8&19.99997634&19.99997634&19.99997636&20.00\\
\bottomrule
\end{tabular}
\end{table}

\section{The case of approximations of distributed order differential operators via asymptotic expansion and GLT momentary symbols}
\label{sec:fractional}

\textcolor{black}{In this last section we focus on the matrix-sequences arising from the numerical approximation of distributed-order operators. In detail, such a procedure consists of two steps:}
\begin{enumerate}
\item  Employ a quadrature formula to discretize the distributed-order operator into a multi-term constant-order fractional derivative;
\item discretize each constant-order fractional derivative.
\end{enumerate}

In particular, we focus on the case where the matrices under consideration take the form
%\begin{equation}\label{eq_sum}
%\mathcal{T}_n:=c_0T_n(f_0)+c_1 h^{\alpha_1} T_n(f_1)+c_2 h^{\alpha_2}T_n(f_2)+\cdots+c_{\tau} h^{\alpha_\tau}T_n(f_{\tau}),
%\end{equation}
\textcolor{black}{
\begin{equation}\label{eq_sum}
\frac{h^{\alpha_\ell}}{\Delta {\alpha}}\mathcal{T}_{n}=c_\ell T_n(g_{\alpha_\ell})+c_{\ell-1} h^{\Delta\alpha} T_n(g_{\alpha_{\ell-1}})+\dots+ c_1h^{\Delta\alpha({\ell-1})}T_n(g_{\alpha_1}),
\end{equation}
where $\ell$ is a positive integer, all the coefficients $c_j$ are positive, independent of $n$, and contained in a specific positive range $[c_*,c^*]$. 
 Moreover, $\Delta\alpha=\frac{1}{\ell}$ and all the terms $1<\alpha_1<\alpha_2\cdots<\alpha_\ell<2$ are positive and defined by $a_k=1+\left(k-\frac{1}{2}\right)\Delta \alpha$, $k=1,\dots, \ell$. More importantly all the functions $g_{\alpha_\ell}$ are globally continuous, monotonically increasing in the interval $[0,\pi]$ and even in the whole definition domain $[-\pi,\pi]$.
 }
 
 \textcolor{black}{
The goal is to exploit the notion of GLT momentary symbols and use it in combination with the asymptotic expansions derived in a quite new research line (see \cite{EkFuSe18} and references there reported), in order to have a very precise description of the spectrum of such matrices.
}

\textcolor{black}{
Indeed, under specific hypotheses on the generating function $f$, and fixing an integer $\nu\ge0$, it is possible to give an accurate description of the eigenvalues of $T_n(f)$ via the following asymptotic expansion:
\begin{align*}
\lambda_j(T_n(f))=w_0(\theta_{j,n})+hw_1(\theta_{j,n})+h^2w_2(\theta_{j,n})+\ldots +h^\nu w_\nu(\theta_{j,n})+E_{j,n,\nu},
\end{align*} 
where  the eigenvalues of $T_n(f)$ are arranged in ascending order,  $h=\frac{1}{n}$, and $\theta_{j,n}=\frac{j\pi}{n}=j\pi h$ for $j=1,\ldots,n$, $E_{j,n,\nu}=O(h^{\nu+1})$ is the error. 
Moreover,  $\{w_k\}_{k=1,2,\ldots}$ is a sequence of functions from $[0,\pi]$ to $\mathbb R$.  The idea of such procedure is that a numerical approximation of the value  $w_k(\theta_{j,n})$ can be obtained by fast interpolation-extrapolation algorithms (see \cite{EkGa2019} and references therein). In particular, choosing  $\nu$ proper grids $\theta_{j,n_1}$, $\theta_{j,n_2}$, 
$\dots$ $\theta_{j,n_\nu}$ with $n>>n_\nu>\dots>n_1$
 an approximation of the quantities $\tilde{w_k}({\theta}_{j,n})\approx w_k({\theta}_{j,n})$ can be obtained.   In the Hermitian case, we find that $\tilde{w}_0$ coincides with the generating function. }
 
 \textcolor{black}{ 
Concerning the example in (\ref{eq_sum}), the idea is to link the  functions $\tilde{w}^{c_i g_{\alpha_i}}_k$, $k=1,\dots,\nu$, associated with each $c_i g_{\alpha_i}$  with the GLT momentary symbols of $\frac{h^{\alpha_\ell}}{\Delta {\alpha}}\mathcal{T}_{n}$. 
Precisely, for $j=1,\dots,n$, for a fixed $\nu$, we approximate the eigenvalues of $\frac{h^{\alpha_\ell}}{\Delta {\alpha}}\mathcal{T}_{n}$ by
}
\textcolor{black}{
\begin{equation}\label{eq:momentary_asymp_exp}
\begin{split}
&\lambda_j\left(\frac{h^{\alpha_\ell}}{\Delta {\alpha}}\mathcal{T}_{n}\right)\approx\\
&c_\ell g_{\alpha_\ell}(\theta_{j,n})+\sum_{t=1}^{\nu}h^t\left(\tilde{w}_t^{{\alpha_\ell}}(\theta_{j,n})+\sum_{i=\ell-1}^1\tilde{w}_{t-1}^{{\alpha_i}}(\theta_{j,n})h^{-(1-\Delta\alpha(\ell-i))}\right),
\end{split}
\end{equation}
where, for sake of notation, we denoted by $\tilde{w}_t^{\alpha_i}$ the approximation of the $t$-th asymptotic expansion coefficient associated with $c_i g_{\alpha_i}$ and the term  $\tilde{w}_0^{\alpha_i}$ coincides with the evaluations of $c_ig_{\alpha_i}$.
 }
 
 \textcolor{black}{
We highlight that the terms in brackets on the right-hand side of the equality act as possible asymptotic expansion coefficients associated with the GLT momentary symbols $g_n$ of $\frac{h^{\alpha_\ell}}{\Delta {\alpha}}\mathcal{T}_{n}$. Note that formula  (\ref{eq:momentary_asymp_exp}) can be rewritten in compact form as
\begin{equation*}
\begin{split}
\lambda_j\left(\frac{h^{\alpha_\ell}}{\Delta {\alpha}}\mathcal{T}_{n}\right)\approx\sum_{t=1}^{\nu}h^t\left(\tilde{w}_t^{{\alpha_\ell}}(\theta_{j,n})+\sum_{i=\ell}^1\tilde{w}_{t-1}^{{\alpha_i}}(\theta_{j,n})h^{-(1-\Delta\alpha(\ell-i))}\right).
\end{split}
\end{equation*}
Hence, it is easy to see that the  GLT momentary symbols correspond to take $\nu$ of the asymptotic expansion equal to 1.
}

\textcolor{black}{
In the following we consider the cases where $\ell=2$ and $\ell=5$ and $\ell=n$ as in Section 4 of \cite{MaSe2021} and confirming at least numerically the conjecture in (\ref{eq:momentary_asymp_exp}) for a fixed $\nu=4$. 
}
\textcolor{black}{
\subsection{Examples:  }
For $\ell=2$, $\Delta\alpha$ is $\frac{1}{2}$ and the matrix in (\ref{eq_sum}) becomes
\[2h^{\alpha_2}\mathcal{T}_{n}=c_2 T_n(g_{\alpha_2})+c_{1} h^{\frac{1}{2}} T_n(g_{\alpha_1}),\]
where $\alpha_1=\frac{5}{4}$ and $\alpha_2=\frac{7}{4}$.
}
\textcolor{black}{
Exploiting the procedure based on formula (\ref{eq:momentary_asymp_exp}) with $\nu=4$, we compute an approximation of the eigenvalues of $c_2 T_n(g_{\alpha_2})+c_{1} h^{\frac{1}{2}} T_n(g_{\alpha_1})$ by 
\begin{equation*}\label{eq:momentary_asymp_exp_l2}
\begin{split}
c_2g_{\alpha_2}({{\theta}_{j,n}})+&h\left[\tilde{w}^{{\alpha_2}}_{1}({\tilde{\theta}_{j,n}})+h^{-\frac{1}{2}}c_1g_{\alpha_1}({{\theta}_{j,n}})\right]+\\
& \sum_{t=2}^{\nu} h^{t}\left[ \tilde{w}^{{\alpha_2}}_{t}({{\theta}_{j,n}})+h^{-\frac{1}{2}}\tilde{w}^{{\alpha_1}}_{t-1}({{\theta}_{j,n}})\right],
\end{split}
\end{equation*}
for $j=1, \dots, n$. We consider the cases where $n=100,500,1000$ using an initial grid with $n_1=10$ points and we compare the aforementioned approximations with those obtained by the evaluations of GLT and GLT momentary symbols associated with the sequence described by the matrices in (\ref{eq_sum}). 
In Figure \ref{fig:frac_expansion_l2} we can observe that the approximation of the spectrum obtained computing the evaluations of the GLT momentary symbols is better with respect to that provided by the evaluations $c_\ell g_{\alpha_\ell}(\theta_{j,n})$. Moreover, the error of the approximation significantly reduces for almost all the eigenvalues  combining the notions of GLT momentary symbols with the asymptotic expansion described before, see Figure \ref{fig:err_frac_expansion_l2}. 
Note that the particular shape of the asymptotic expansion error depends on the fact that in correspondence with the grid points $\theta_{j,n_t}$, $t=1,\dots, \nu$ the quantities $\tilde{w}_t^{\alpha_i}$ are calculated exactly by the extrapolation-interpolation procedure. Moreover, note that the accuracy of the approximation via the combination of GLT momentary symbols and spectral asymptotic expansion seems to decrease corresponding to the maximum eigenvalue. Actually, this behavior is expected by the theory of the asymptotic expansion. Indeed, it is a consequence of the fact that the involved symbols are not trigonometric polynomials and in particular they become non-smooth when periodically extended on the real line.
}

\begin{figure}[!h]
\centering
\includegraphics[width=0.49\textwidth]{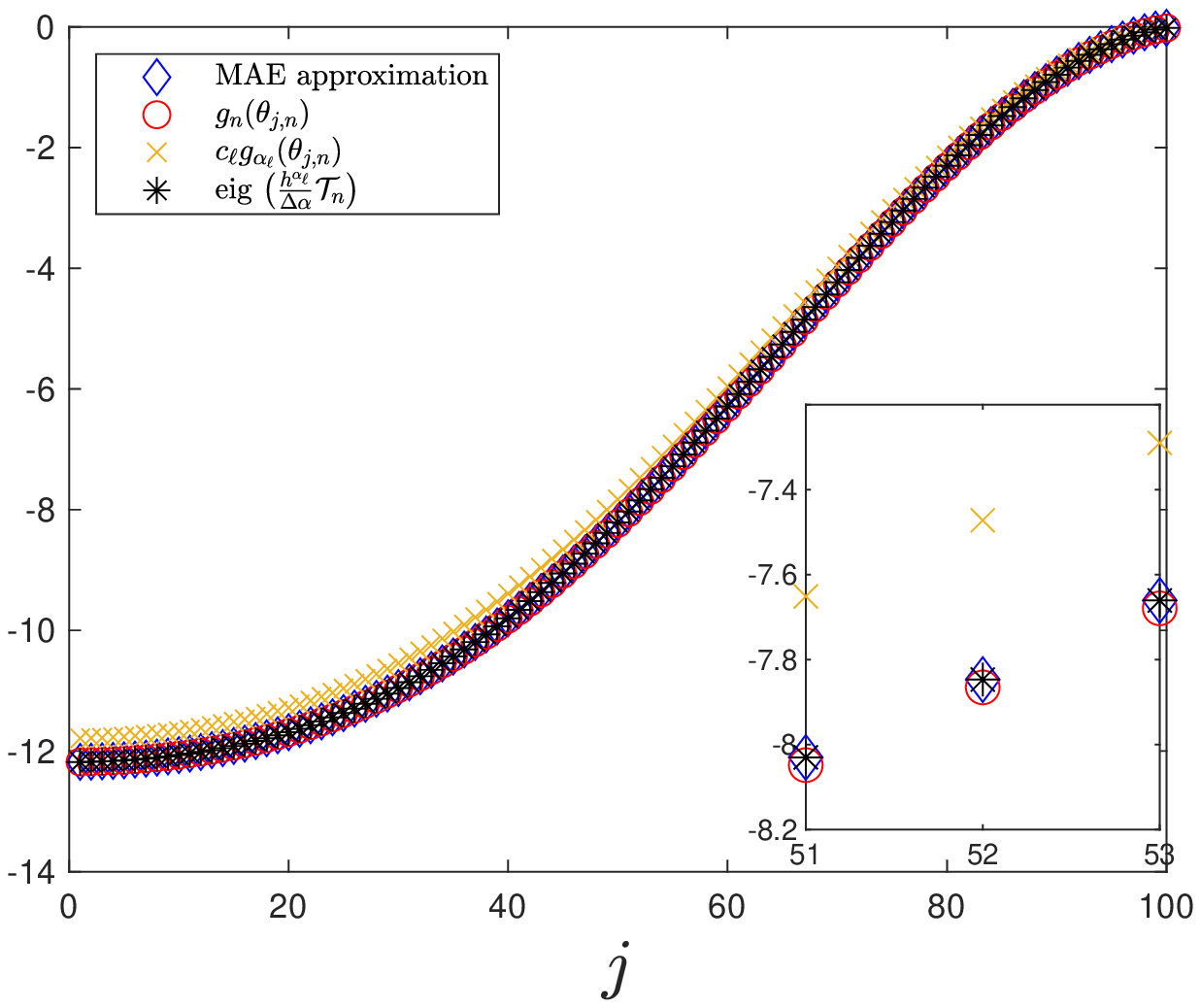}\\
\includegraphics[width=0.49\textwidth]{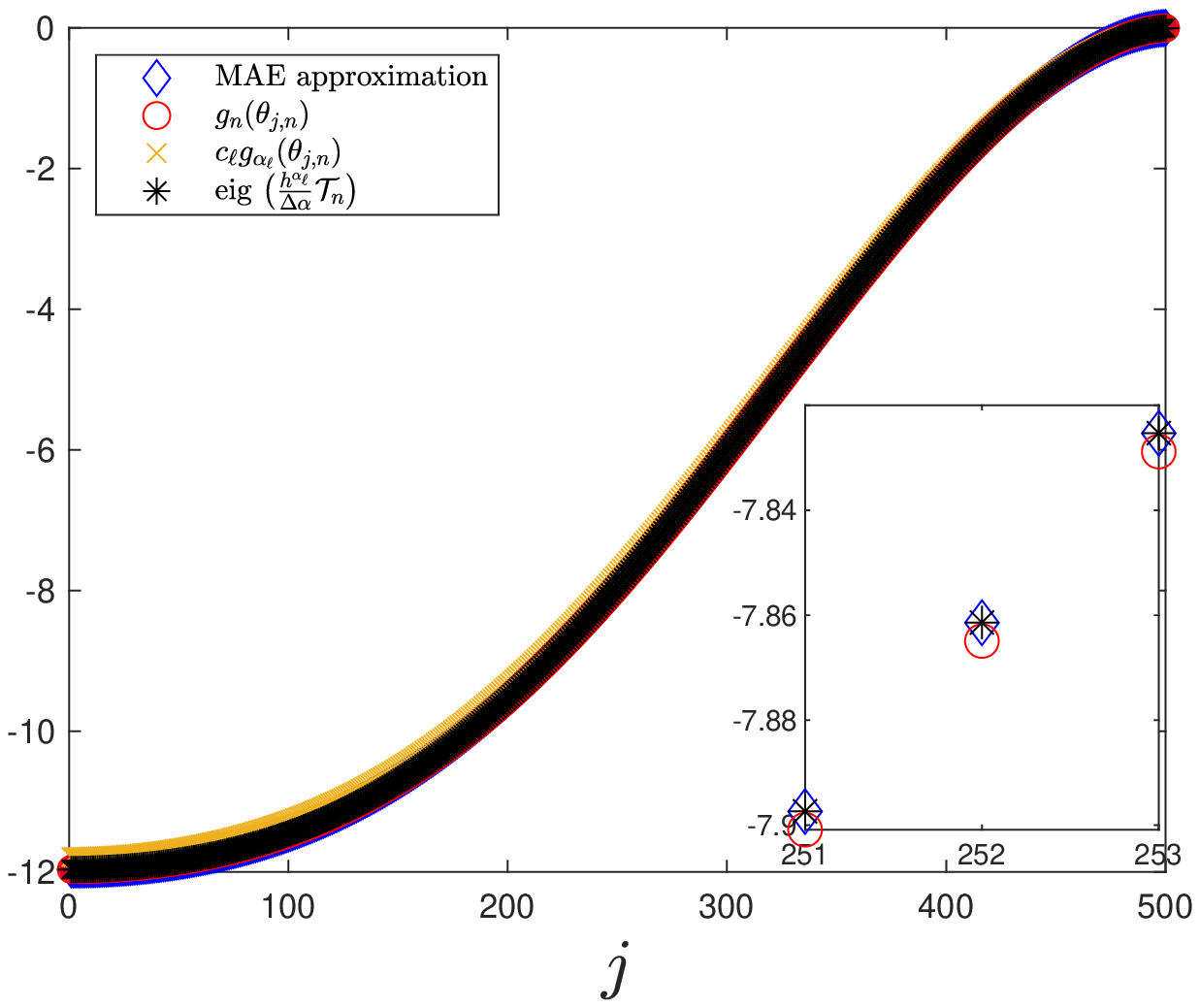}
\includegraphics[width=0.49\textwidth]{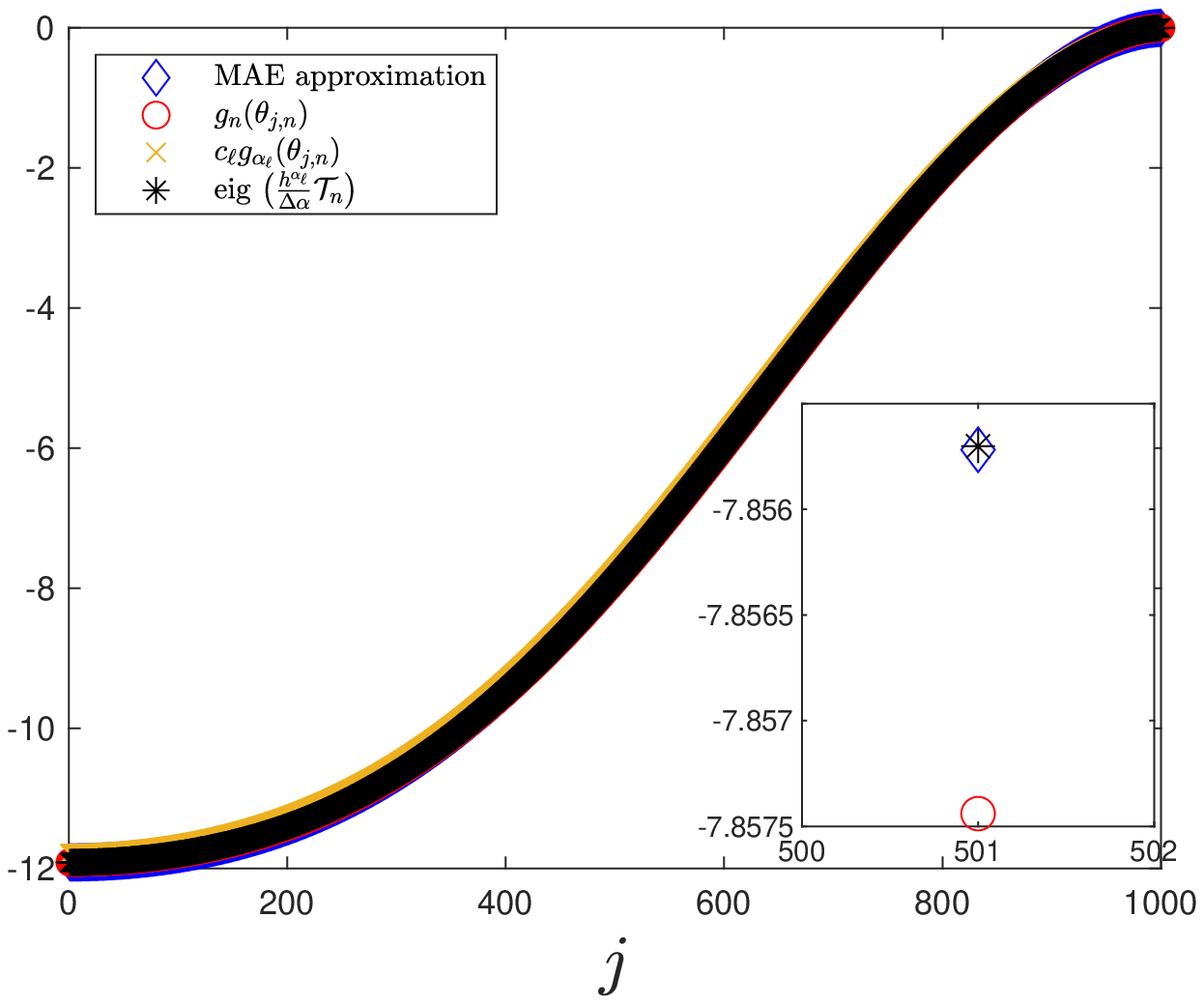}
\caption{\textcolor{black}{Approximation of the eigenvalues of  $\frac{h^{\alpha_\ell}}{\Delta \alpha}\mathcal{T}_{n}$, $\ell=2$ by the samplings of the GLT and GLT momentary symbols, and making use of the momentary asymptotic expansion (MAE) with $\nu=4$ for $n=100,500,1000$ with an initial grid of $n_1=10$ points. }}
\label{fig:frac_expansion_l2}
\end{figure}

\begin{figure}[!h]
\centering
\includegraphics[width=0.6\textwidth]{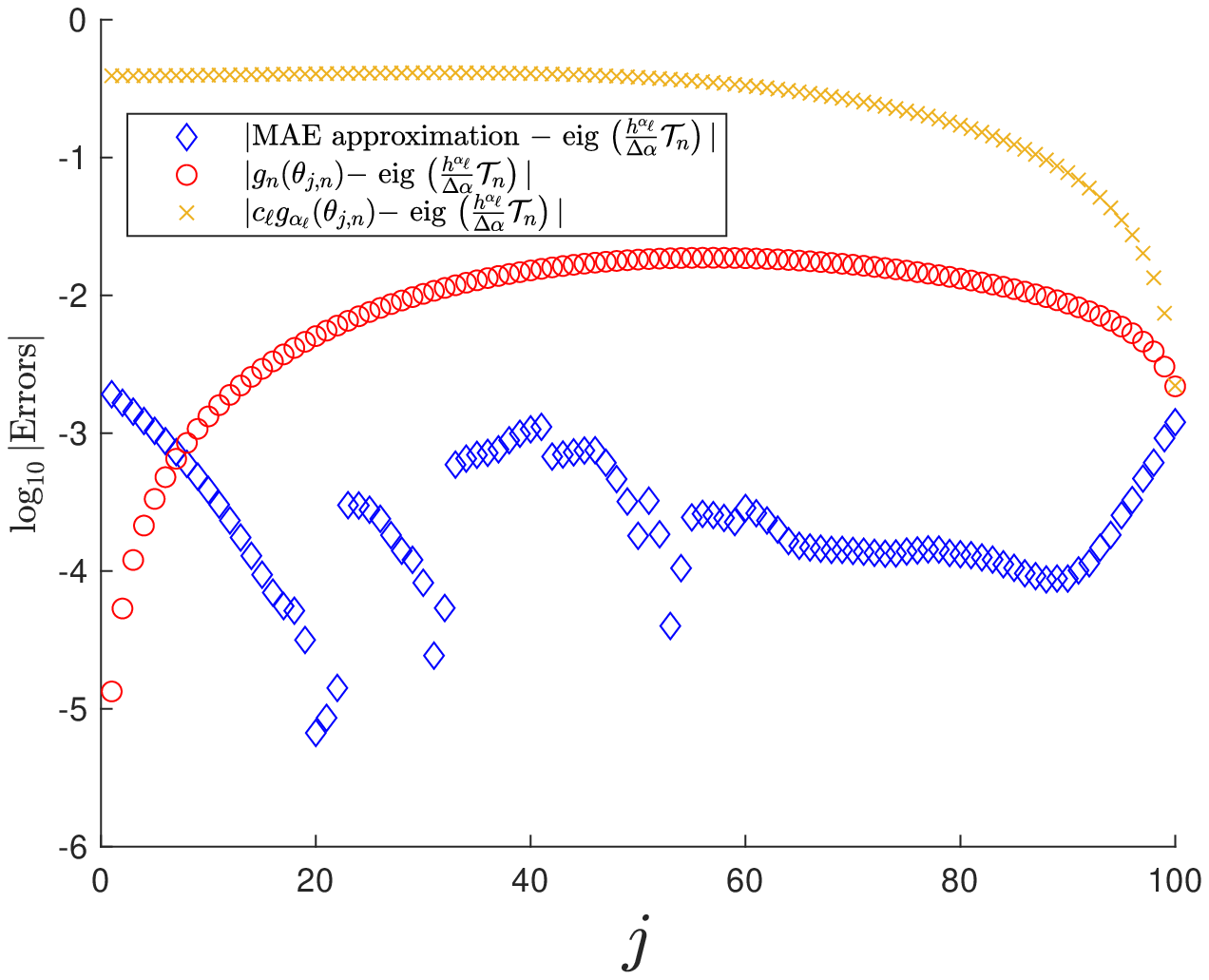}\\
\includegraphics[width=0.49\textwidth]{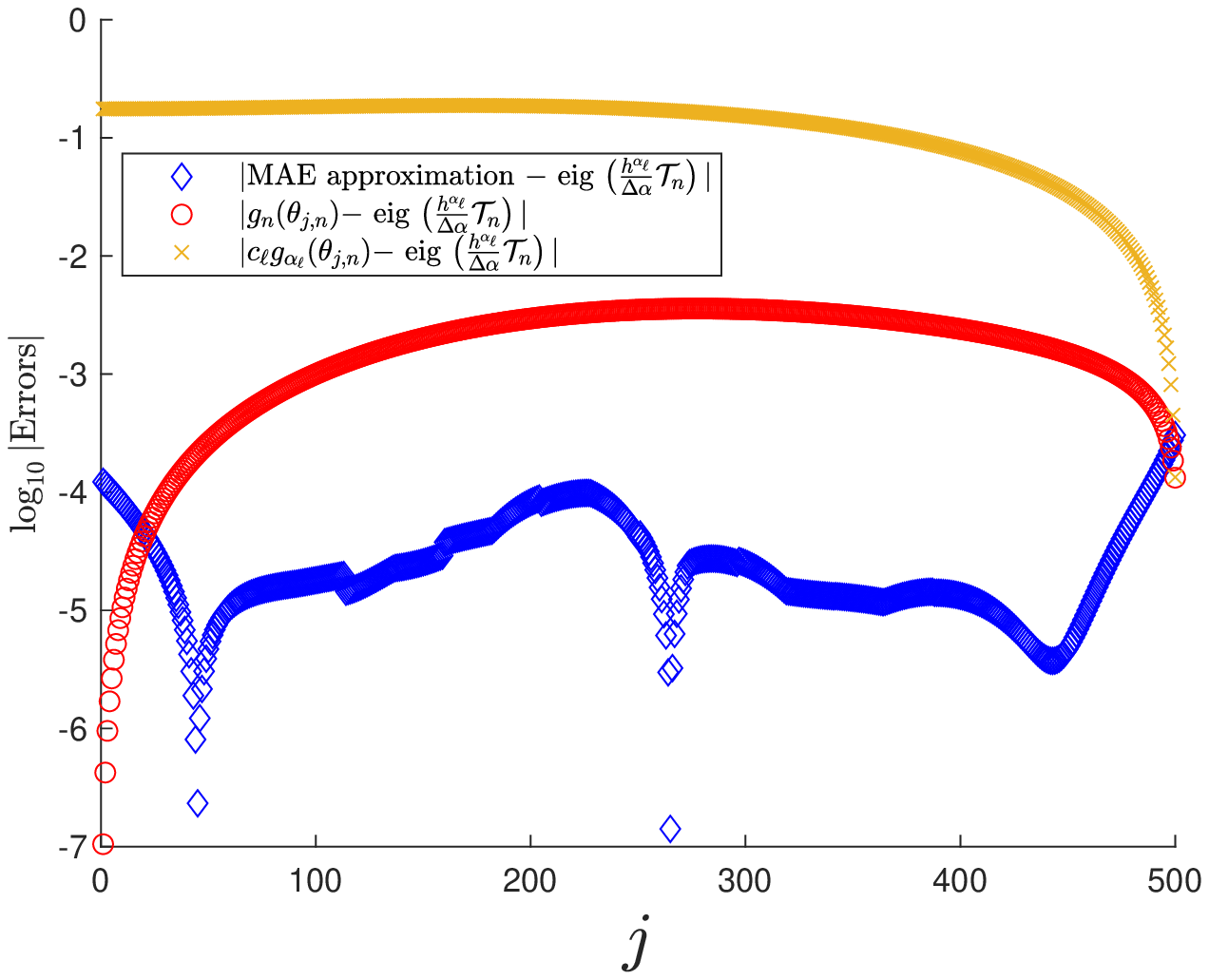}
\includegraphics[width=0.49\textwidth]{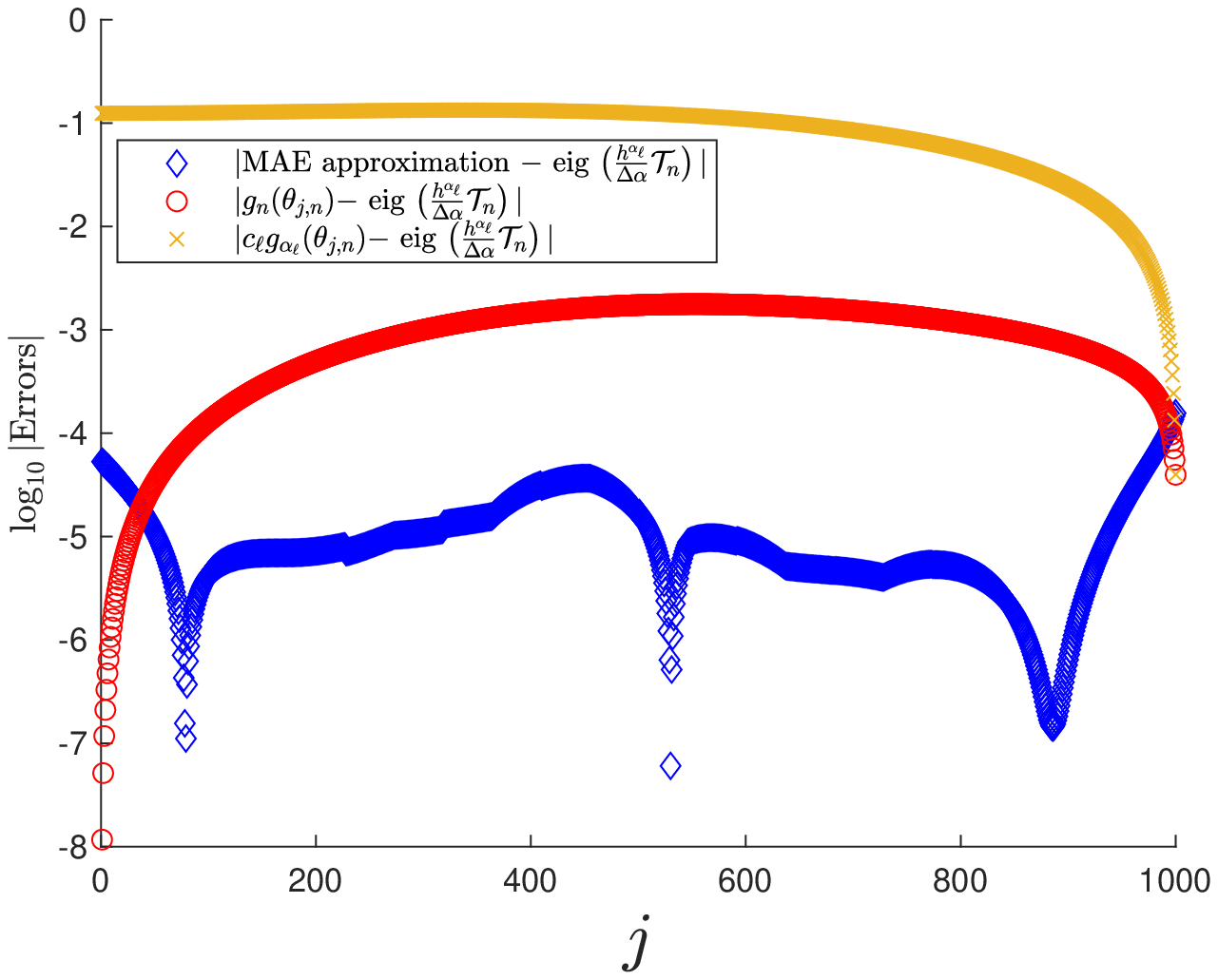}
\caption{\textcolor{black}{Absolute errors of the approximation of the eigenvalues of  $\frac{h^{\alpha_\ell}}{\Delta \alpha}\mathcal{T}_{n}$, $\ell=2$ by the samplings of the GLT and GLT momentary symbols, and making use of the  momentary asymptotic expansion (MAE) with $\nu=4$ for $n=100,500,1000$ with an initial grid of $n_1=10$ points. }}
\label{fig:err_frac_expansion_l2}
\end{figure}

\textcolor{black}{
Following the analogous procedure, we consider the case where $\ell=5$ and $\ell=n$ which are associated with $\Delta \alpha=\frac{1}{5}$ and  $\Delta \alpha=\frac{1}{n}$, respectively. In Figures \ref{fig:frac_expansion_l5} and  \ref{fig:frac_expansion_ln} we plot the approximations of the eigenvalues given by the three presented strategies for $\ell=5$ and $\ell=n$. }
\textcolor{black}{
Again, we obtain numerical confirmation that the combination of the notions of GLT momentary symbols and asymptotic expansion provides accurate results even for moderate sizes, as confirmed by the error plots in Figures \ref{fig:err_frac_expansion_l5} and  \ref{fig:err_frac_expansion_ln}, for $n=100,500,1000$.
}
The good outcome of the presented numerical tests gives ground for a finer analysis of the spectral features of the matrices considered in the case  left open in \cite{MaSe2021}, which arises when the integral partition width is asymptotic to the adopted discretization step. That is, when in formula (\ref{eq_sum}) we take $\alpha_j=jh$, $j=1,\ldots,\ell$, \textcolor{black}{$\ell=n$}. 

Moreover, efficient and fast algorithms  which exploit the concept of momentary symbols can be studied for computing the singular values and eigenvalues of $T_n(f)$ with its possible block, and variable coefficients generalizations and this will be investigated in the future.

\begin{figure}[!h]
\centering
\includegraphics[width=0.49\textwidth]{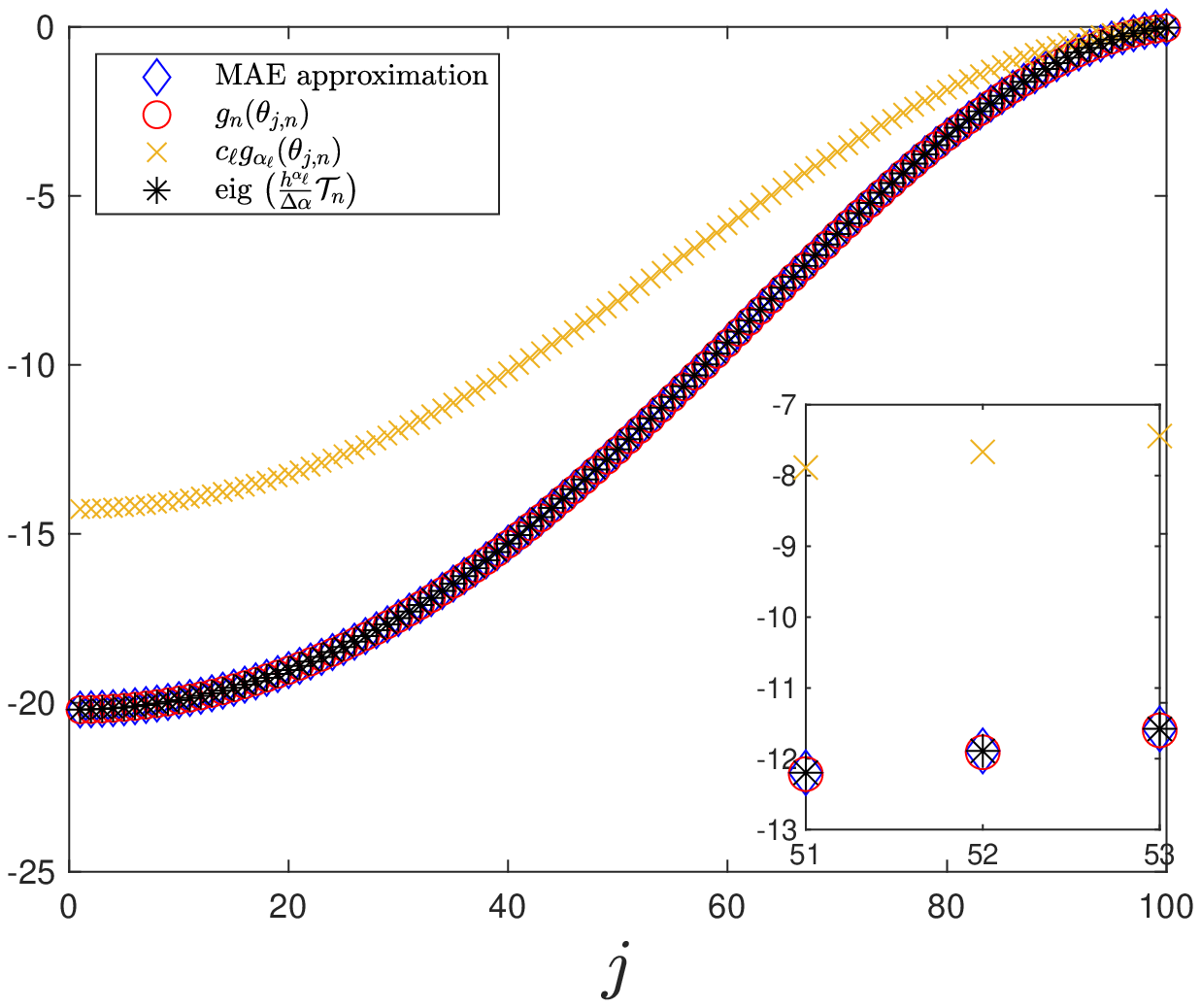}\\
\includegraphics[width=0.49\textwidth]{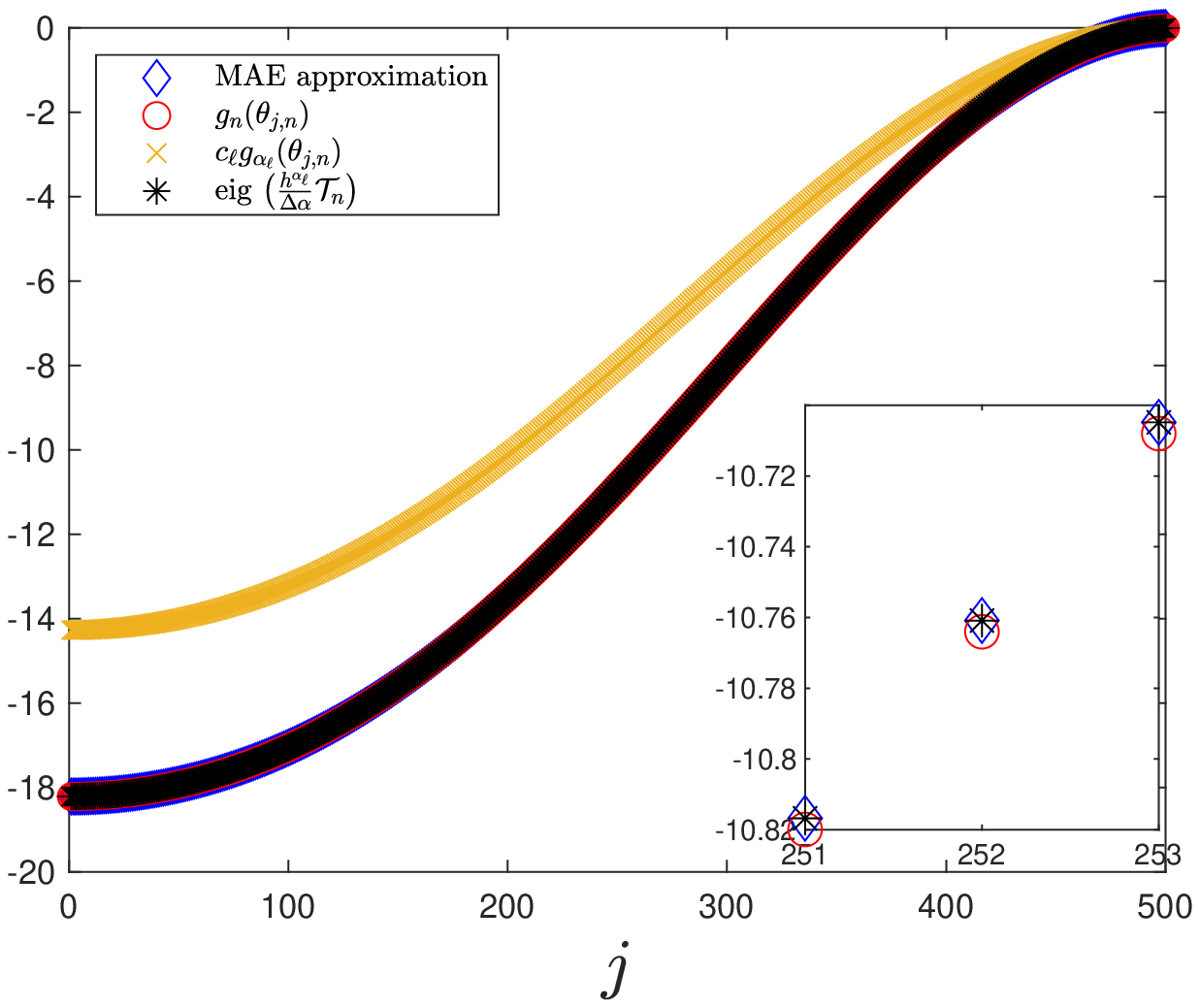}
\includegraphics[width=0.49\textwidth]{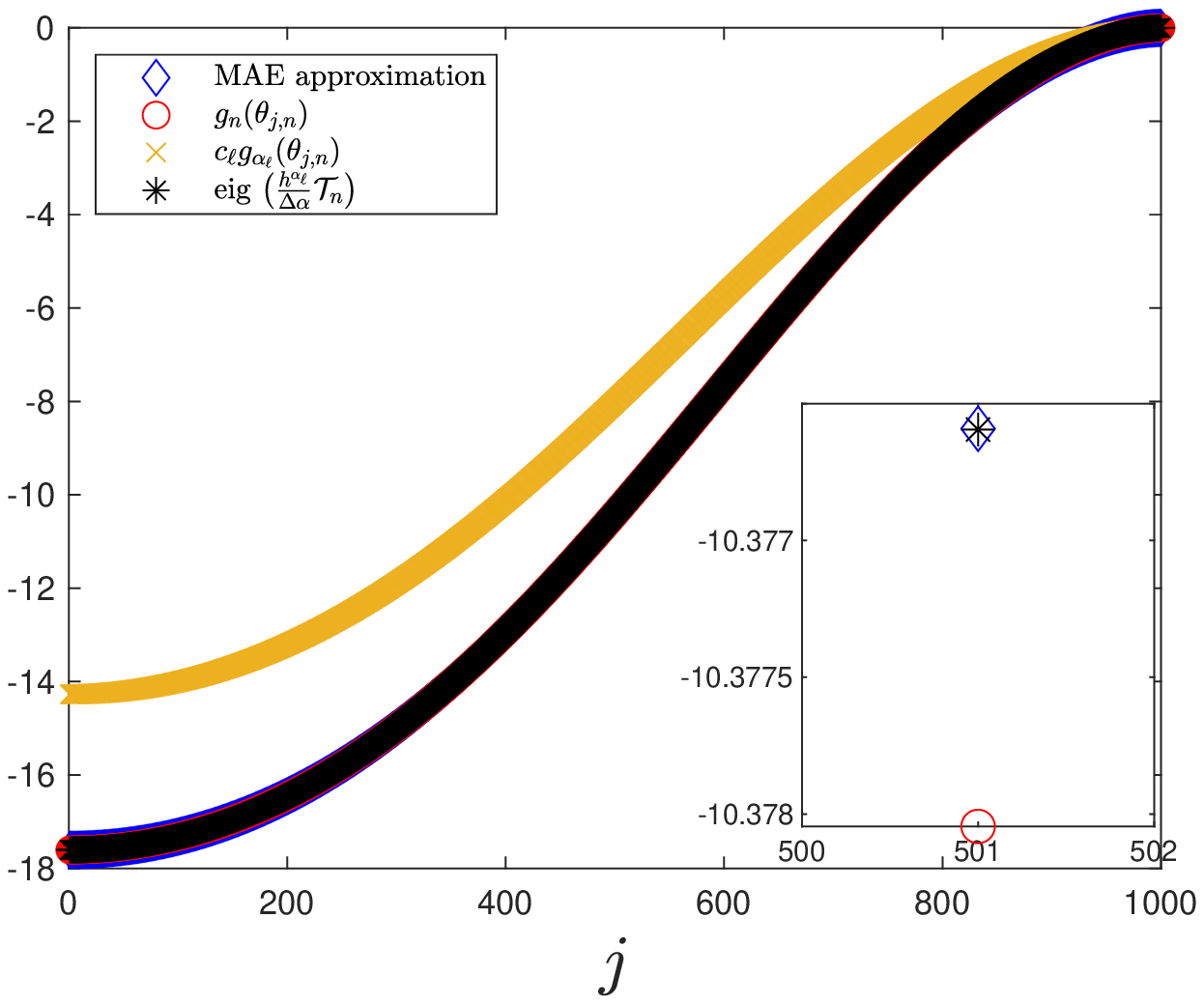}
\caption{\textcolor{black}{Approximation of the eigenvalues of  $\frac{h^{\alpha_\ell}}{\Delta \alpha}\mathcal{T}_{n}$, $\ell=5$ by the samplings of the GLT and GLT momentary symbols, and making use of the  momentary asymptotic expansion (MAE) with $\nu=4$ for $n=100,500,1000$ with an initial grid of $n_1=10$ points. }}
\label{fig:frac_expansion_l5}
\end{figure}
\begin{figure}[!h]
\centering
\includegraphics[width=0.49\textwidth]{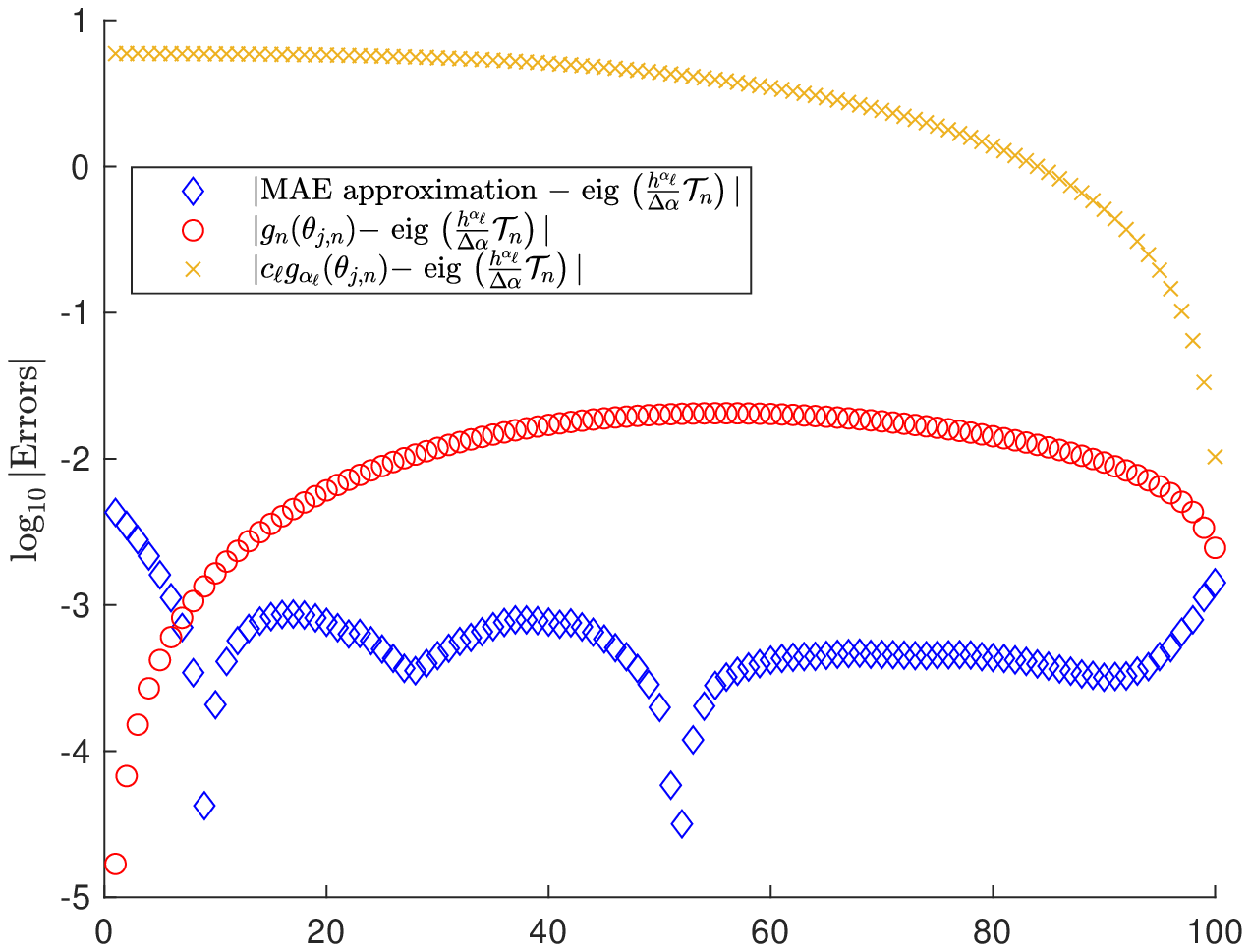}\\
\includegraphics[width=0.49\textwidth]{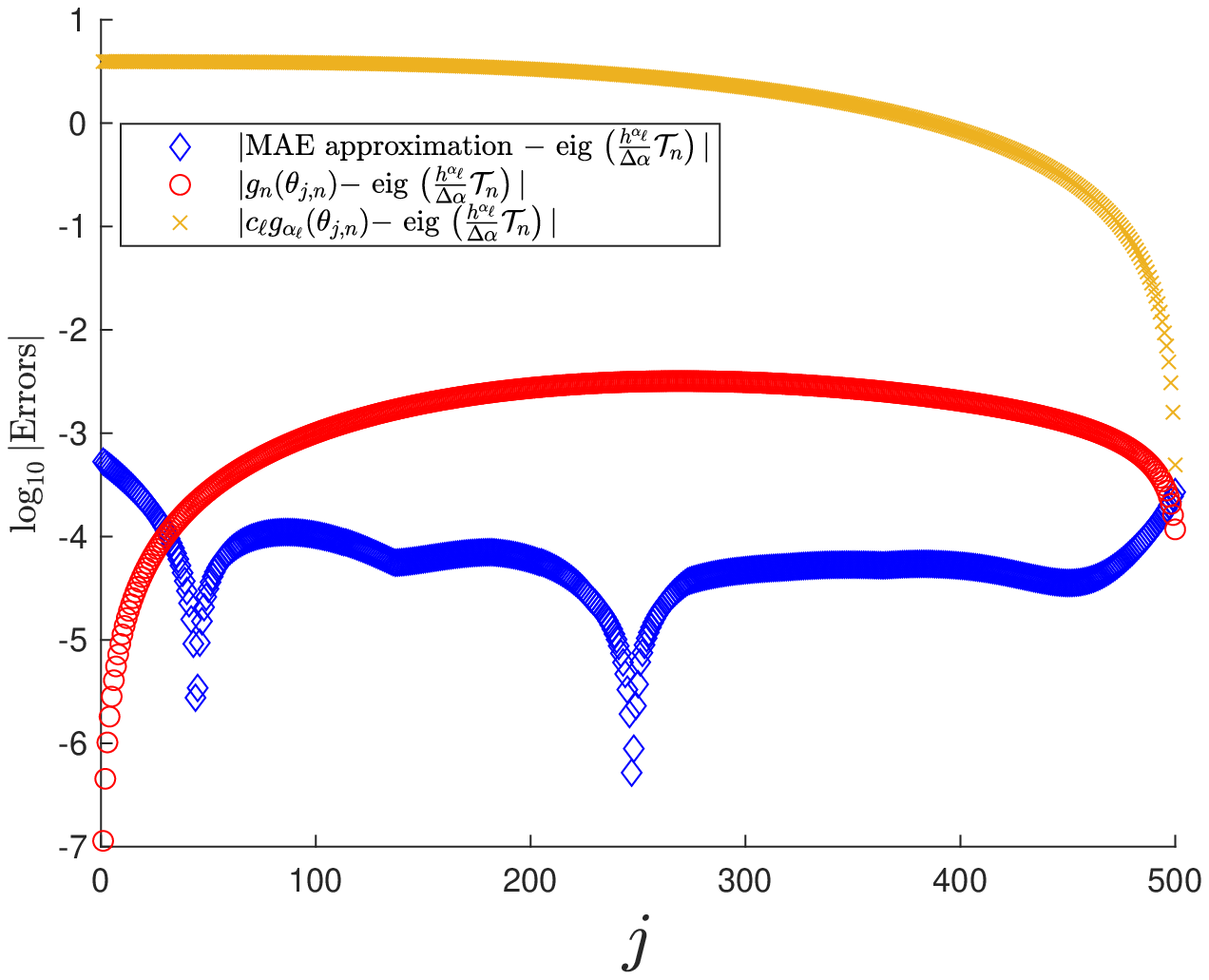}
\includegraphics[width=0.49\textwidth]{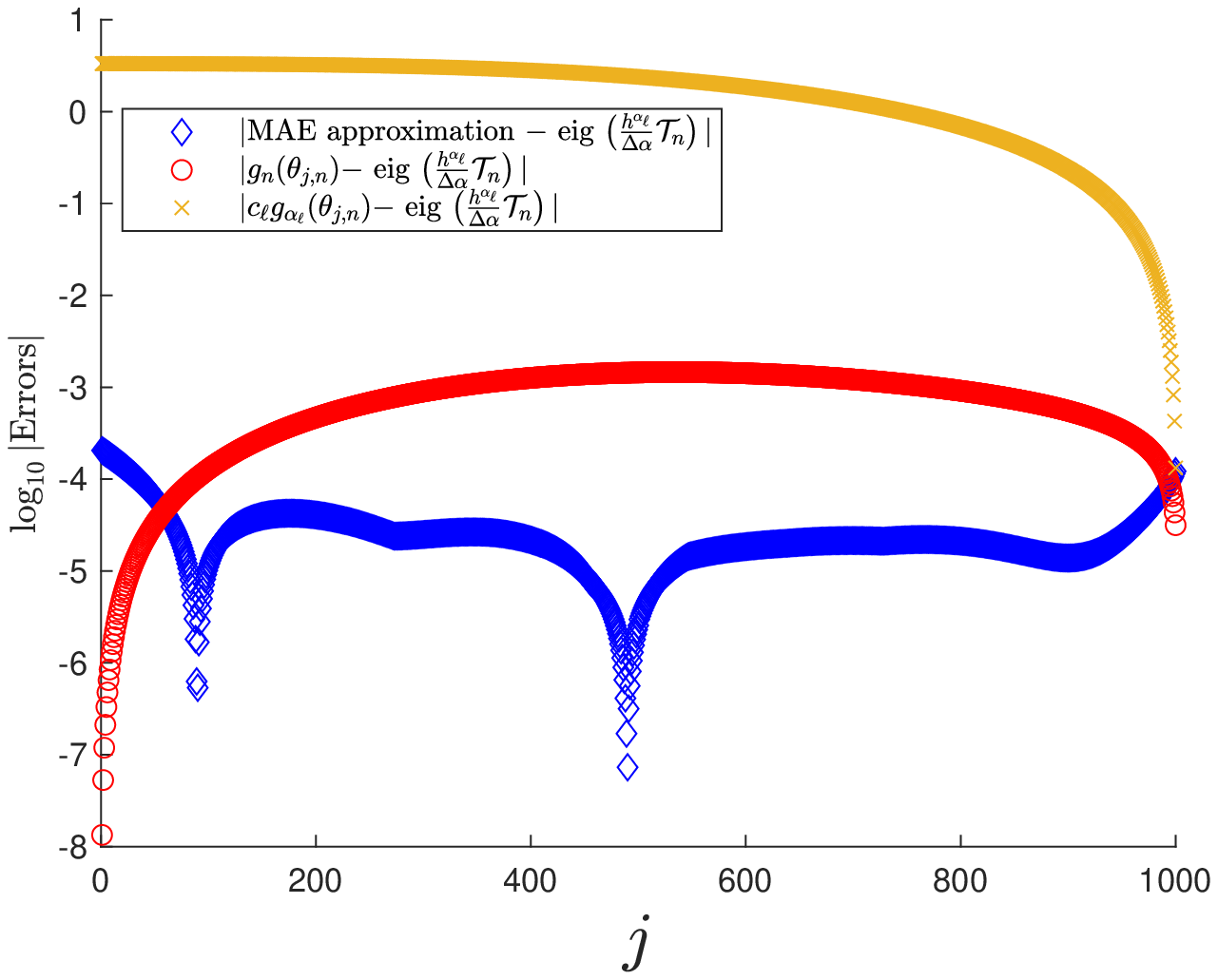}
\caption{\textcolor{black}{Absolute errors of the approximation of the eigenvalues of  $\frac{h^{\alpha_\ell}}{\Delta \alpha}\mathcal{T}_{n}$, $\ell=5$ by the samplings of the GLT and GLT momentary symbols, and making use of the  momentary asymptotic expansion (MAE) with $\nu=4$ for $n=100,500,1000$ with an initial grid of $n_1=10$ points. }}
\label{fig:err_frac_expansion_l5}
\end{figure}

\begin{figure}[!h]
\centering
\includegraphics[width=0.49\textwidth]{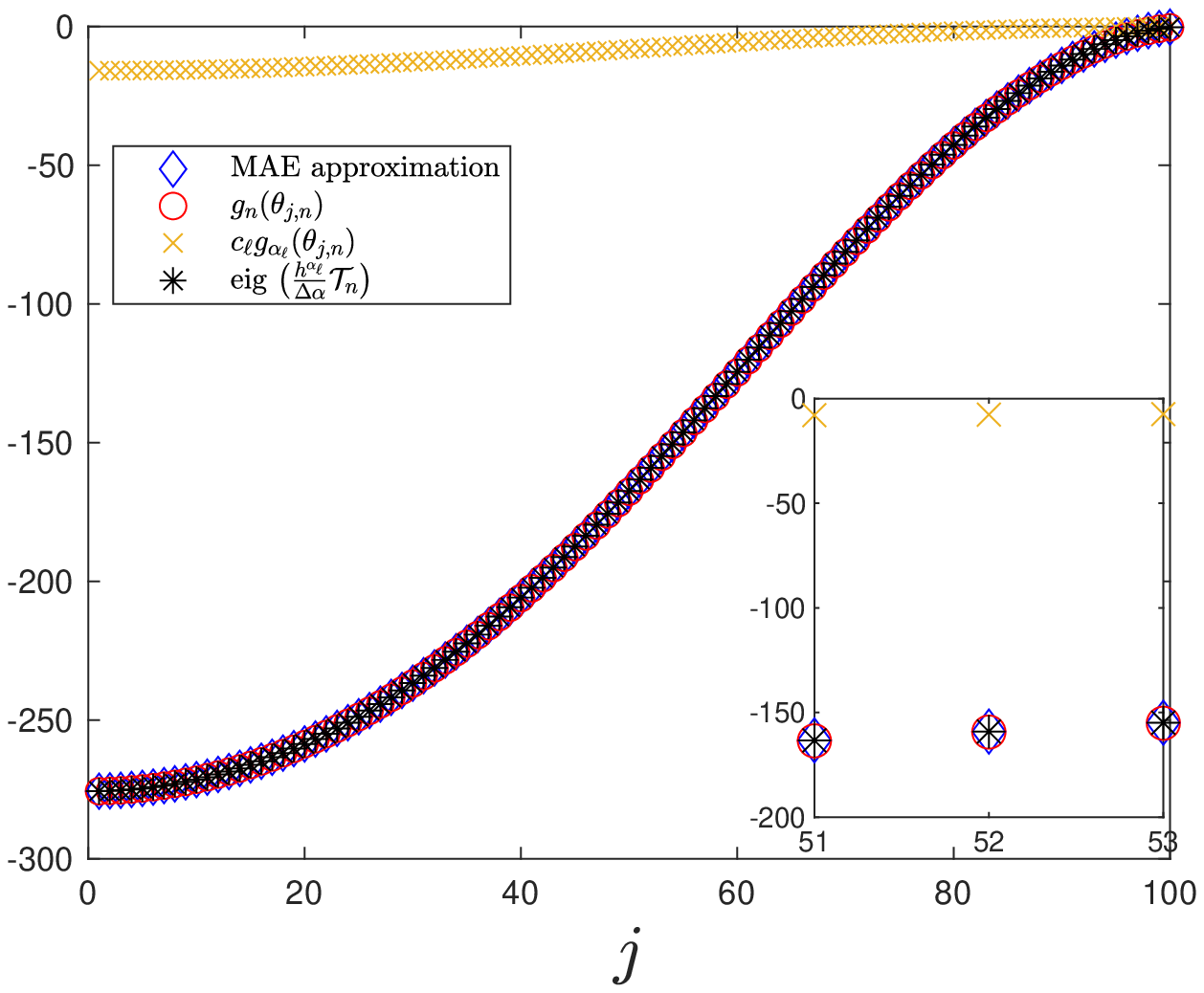}\\
\includegraphics[width=0.49\textwidth]{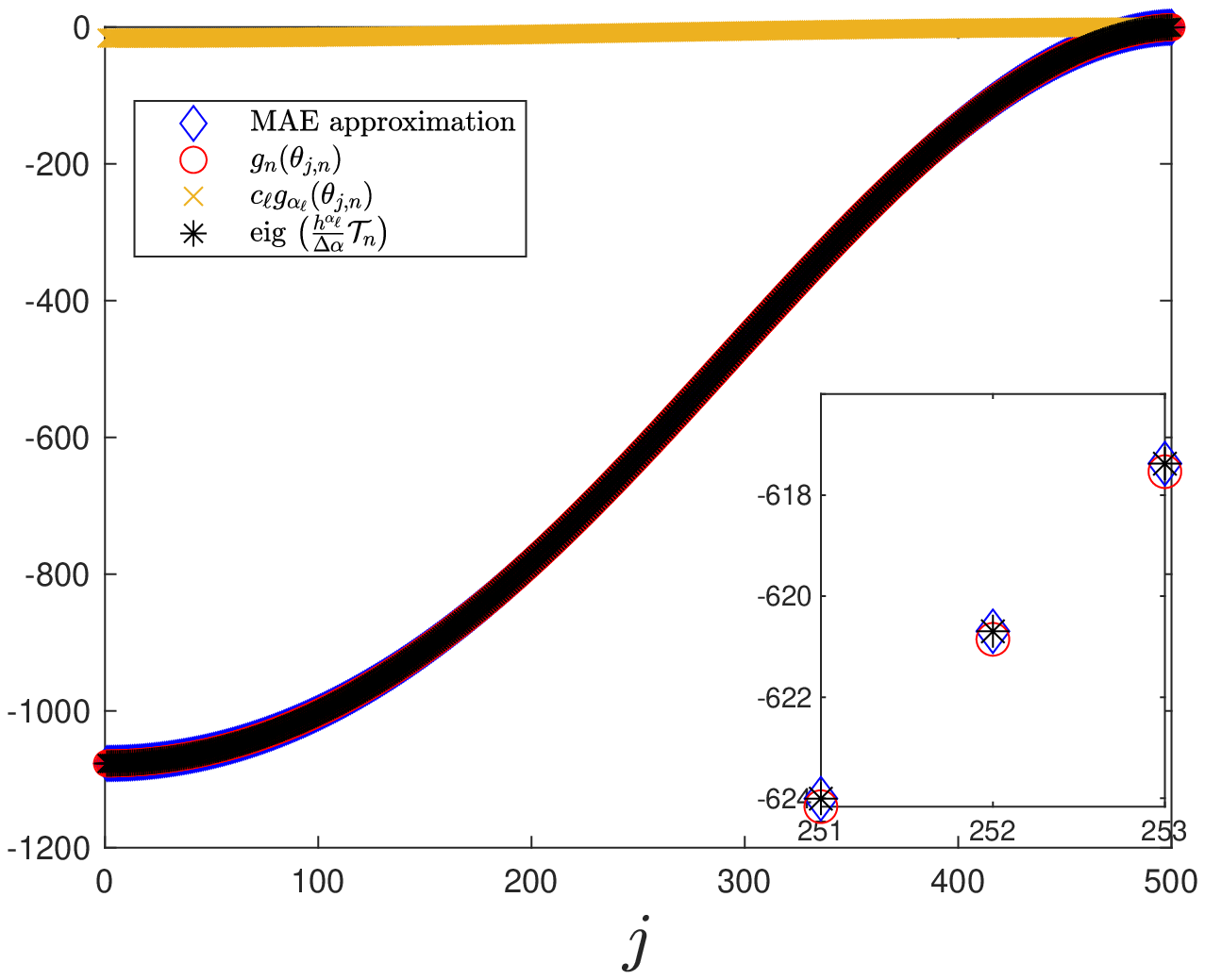}
\includegraphics[width=0.49\textwidth]{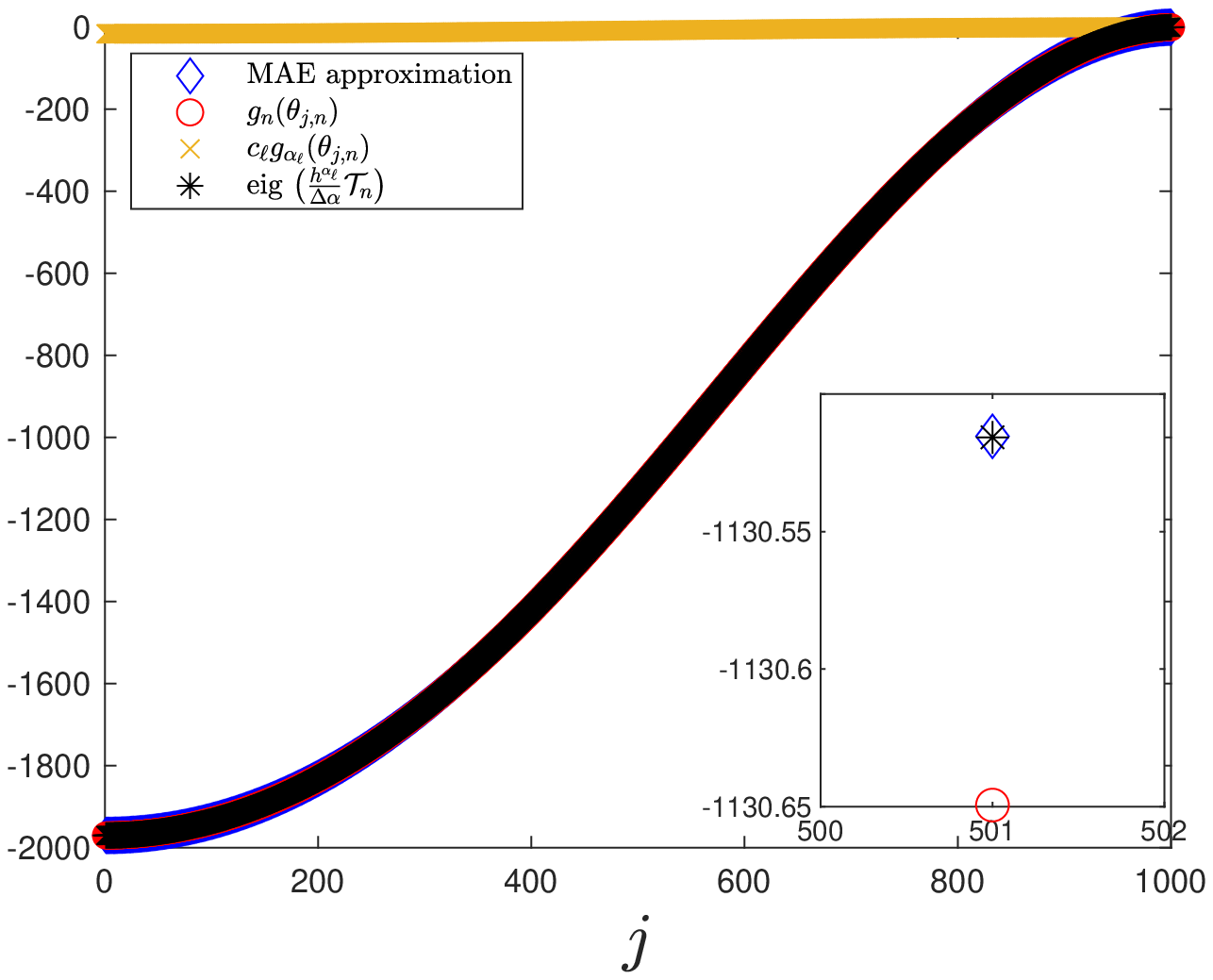}
\caption{\textcolor{black}{Approximation of the eigenvalues of  $\frac{h^{\alpha_\ell}}{\Delta \alpha}\mathcal{T}_{n}$, $\ell=n$ by the samplings of the GLT and GLT momentary symbols, and making use of the  momentary asymptotic expansion (MAE) with $\nu=4$ for $n=100,500,1000$ with an initial grid of $n_1=10$ points. }}
\label{fig:frac_expansion_ln}
\end{figure}

\begin{figure}[!h]
\centering
\includegraphics[width=0.49\textwidth]{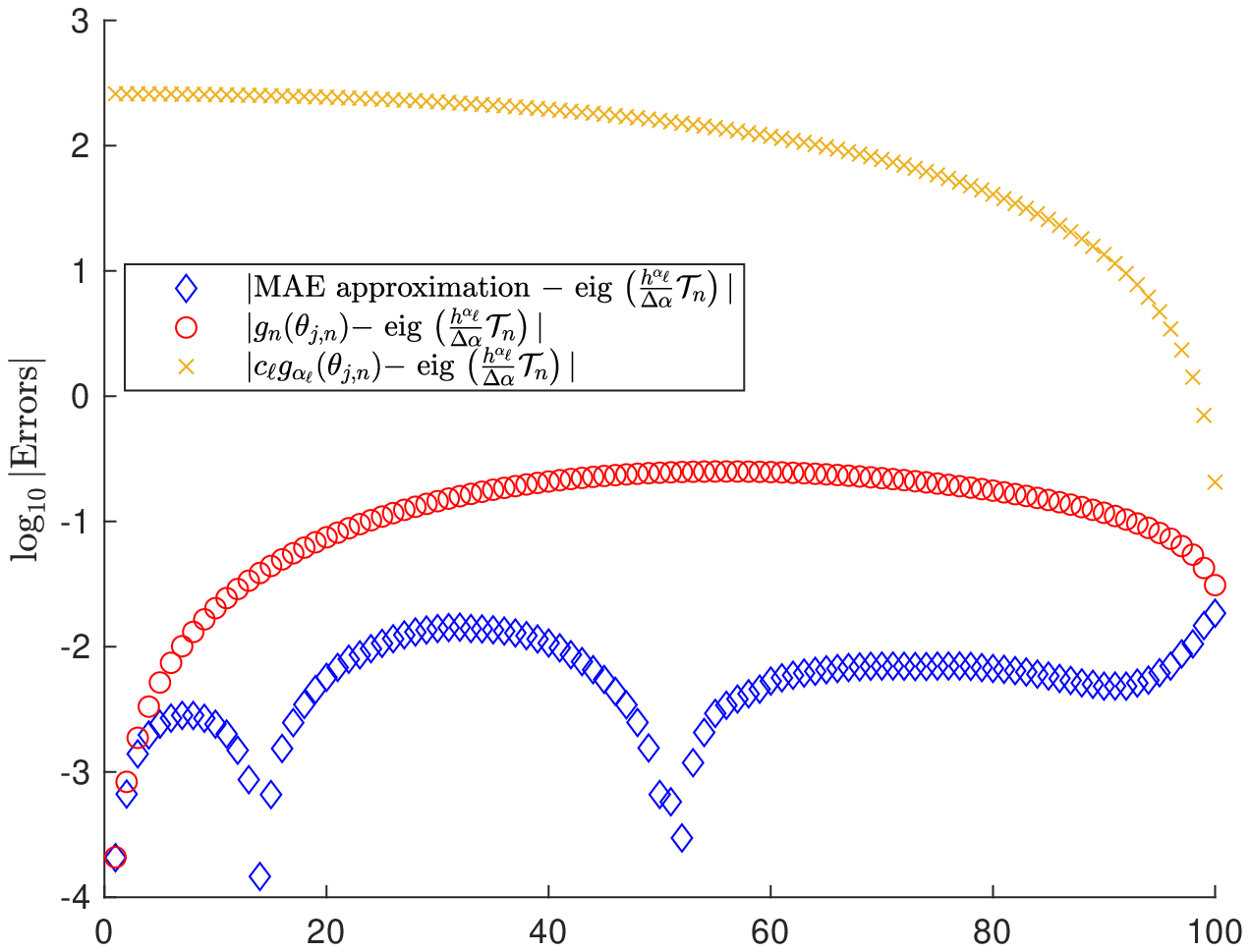}\\
\includegraphics[width=0.49\textwidth]{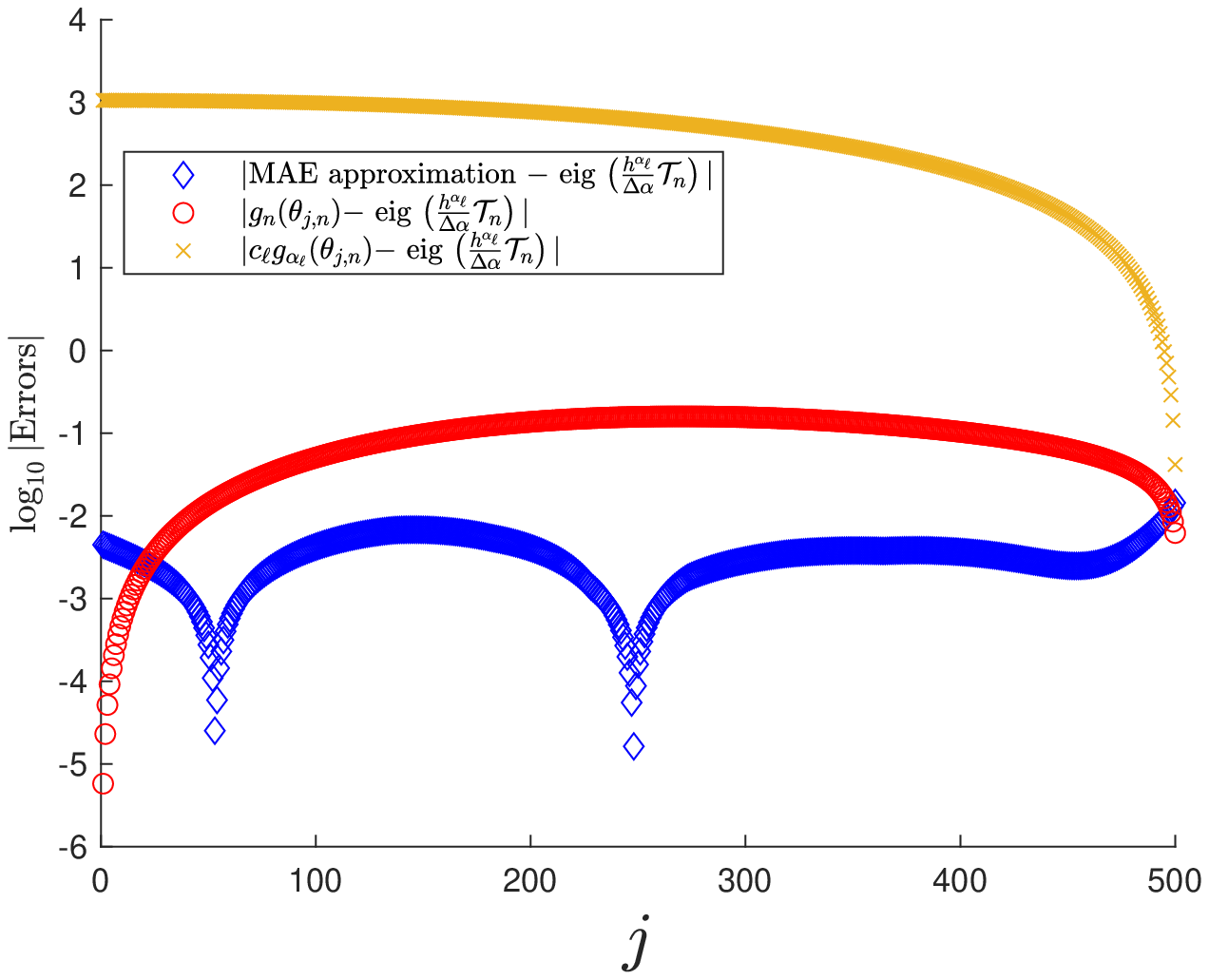}
\includegraphics[width=0.49\textwidth]{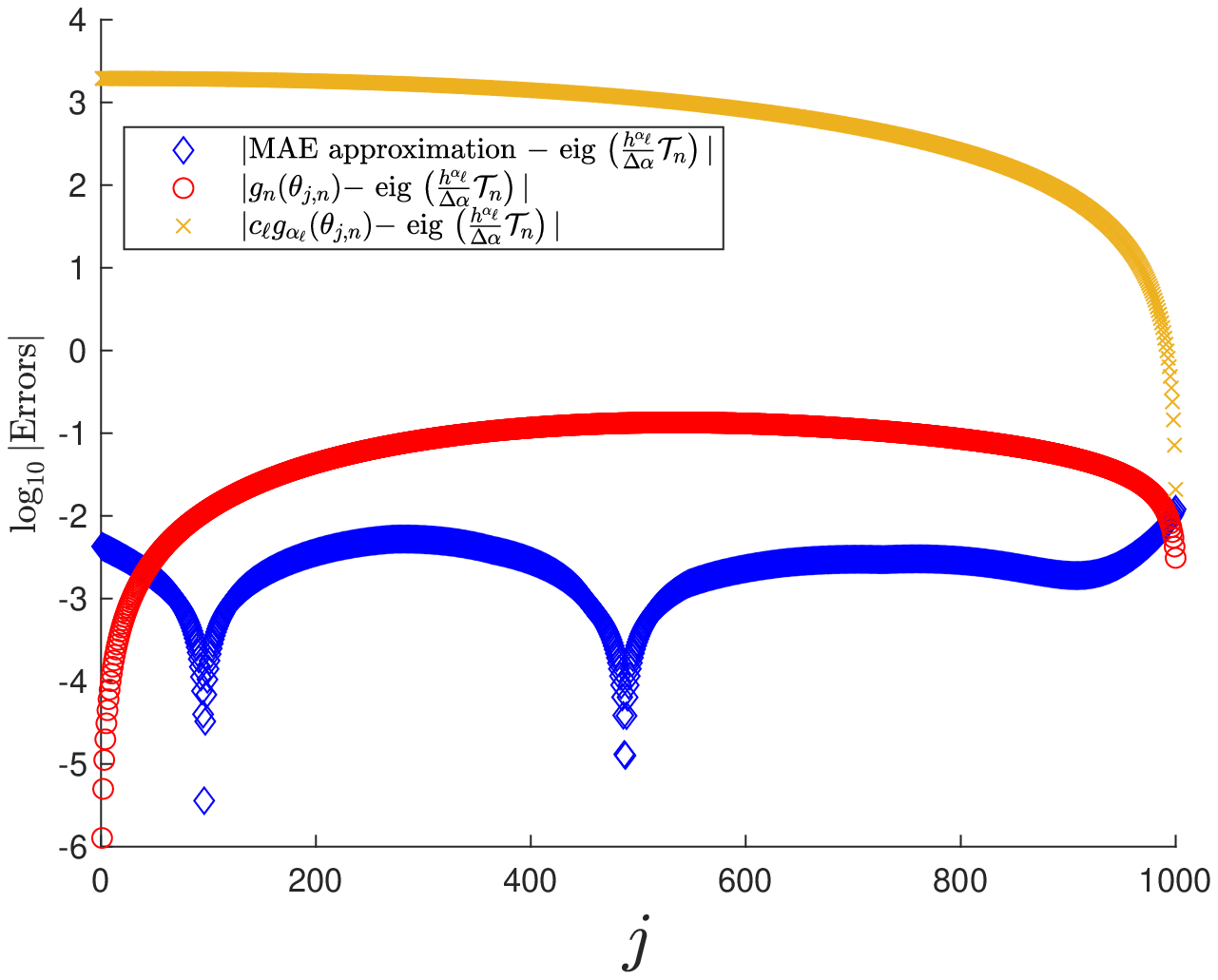}
\caption{\textcolor{black}{Absolute errors of the approximation of the eigenvalues of  $\frac{h^{\alpha_\ell}}{\Delta \alpha}\mathcal{T}_{n}$, $\ell=n$ by the samplings of the GLT and GLT momentary symbols, and making use of the   momentary asymptotic expansion (MAE) with $\nu=4$ for $n=100,500,1000$ with an initial grid of $n_1=10$ points. }}
\label{fig:err_frac_expansion_ln}
\end{figure}

\section{Concluding remarks}
\label{sec:conclusions}

The main focus of this paper has been the characterization of the spectrum and the singular values of the coefficient matrix stemming from the approximation with space-time grid for a parabolic diffusion problem and from the approximation of distributed order fractional equations.
For this purpose we employed the classical GLT theory and the new concept of GLT momentary symbols. The first has permitted to describe the singular value or eigenvalue asymptotic distribution of the sequence of the coefficient matrices. The latter has permitted to derive a function, able to describe the singular value or eigenvalue distribution of the matrix of the sequence, even for small matrix-sizes, but under given assumptions. \textcolor{black}{In particular, we exploited the notion of GLT momentary symbols and we used it in combination with the interpolation-extrapolation algorithms based on the spectral asymptotic expansion of the involved matrices.} 

Many questions remain and below we list open problems to be considered in future researches.
\begin{itemize}
\item  More examples of the use of GLT momentary symbols in a non-Toeplitz setting;
\item The application of GLT momentary symbol in a pure Toeplitz setting, but of very involved nature, like that expressed in relation (\ref{eq_sum}). The use of GLT momentary symbol for the analysis of efficient iterative solvers, also of multigrid type, of linear systems as those appearing in (\ref{eq:system}), also with the inclusion of variable coefficients.
\end{itemize}

\subsection*{Acknowledgment}
This work  was partially supported by INdAM-GNCS. Moreover, the work of Isabella Furci was also supported by the Young Investigator Training Program 2020 (YITP 2019) promoted by ACRI.

% ------------------------------------------------------------------------
\end{document}